\documentclass[12pt]{article}

\usepackage{amsmath,amssymb,amsthm}
\usepackage{graphicx}
\usepackage{hyperref}

\title{The Symmetric Subset Problem in Continuous Ramsey Theory}
\author {Greg Martin\thanks{University of British Columbia, \url{gerg@math.ubc.ca},
        supported in part by grants from the Natural Sciences and Engineering Research Council.}\
        \ and\
        {Kevin O'Bryant\thanks{The City University of New York, College of Staten Island, \url{kevin@member.ams.org},
        supported by NSF grant DMS-0202460.}}
        }
\date{\today}

\newcommand{\Z}{{\mathbb Z}}

\newcommand{\R}{{\mathbb R}}

\newcommand{\K}{{\mathbb K}}
\newcommand{\T}{{\mathbb T}}

\newcommand{\A}{A}
\newcommand{\Prob}{\mathop{\rm Pr}}

\newcommand{\vecd}{{\bf d}}
\newcommand{\vecz}{{\bf z}}
\newcommand{\e}{\varepsilon}
\newcommand{\De}{\Delta(\varepsilon)}

\newcommand{\ff}{f\ast f}
\newcommand{\Bg}{\mbox{$B^{\ast}[g]$} }
\newcommand{\Bgm}{\mbox{$B^{\ast}[g]\pmod{n}$} }
\newcommand{\ffi}{\|\ff\|_\infty}

\newcommand{\Dconstant}{0.591389}
\newcommand{\twotimesDconstant}{1.182778} 
\newcommand{\Rconstant}{1.30036} 
\newcommand{\TrivialLowerBoundconstant}{0.61522} 
\newcommand{\fstarftwonormconstant}{1.14915} 

\newcommand{\floor}[1]{\left\lfloor #1 \right\rfloor}

\newcommand{\lnorm}[3]{{}_{#1}\| #3 \|_{#2}}

\newcommand{\sdr}[1]{ {#1}^{\text{sdr}}}

\newtheorem{thm}{Theorem}[section]
\newtheorem{lem}[thm]{Lemma}

\newtheorem{cor}[thm]{Corollary}

\newtheorem{cnj}[thm]{Conjecture}
\newtheorem{prop}[thm]{Proposition}

\begin{document}
    \maketitle

\begin{abstract}
A symmetric subset of the reals is one that remains invariant under
some reflection $x\mapsto c-x$. We consider, for any $0<\e \le 1$,
the largest real number $\De$ such that every subset of $[0,1]$ with
measure greater than $\e$ contains a symmetric subset with measure
$\De$. In this paper we establish upper and lower bounds for $\De$
of the same order of magnitude: for example, we prove that
$\De=2\e-1$ for $\frac{11}{16}\le\e\le1$ and that
$0.59\e^2<\De<0.8\e^2$ for $0<\e\le\frac{11}{16}$.

This continuous problem is intimately connected with a corresponding discrete problem. A set $S$ of integers is
called a $\Bg$ set if for any given $m$ there are at most $g$ ordered pairs $(s_1,s_2)\in S \times S$ with
$s_1+s_2=m$; in the case $g=2$, these are better known as Sidon sets. Our lower bound on $\De$ implies that
every $\Bg$ set contained in $\{1,2,\dots,n\}$ has cardinality less than $1.30036\sqrt{gn}$. This improves a
result of Green for $g\ge 30$. Conversely, we use a probabilistic construction of $\Bg$ sets to establish an
upper bound on $\De$ for small $\e$.
\end{abstract}

\noindent {\bf AMS Msc (2000):}
    \begin{itemize}
    \item 05D99 Extremal Combinatorics,
    \item 42A16 Fourier Series of Functions with special properties,
    \item 11B83 Special Sequences.
    \end{itemize}


\thispagestyle{empty}   \tableofcontents

\section{Introduction}

A set $C\subseteq\R$ is {\em symmetric} if there exists a number $c$ (the center of $C$)
such that $c+x\in C$ if and only if $c-x\in C$. Given a set $A\subseteq[0,1)$ of
positive measure, is there necessarily a symmetric subset $C\subseteq A$ of positive
measure? The main topic of this paper is to
determine how large, in terms of the Lebesgue measure of $A$, one may take the symmetric
set $C$. In other words, for each $\e>0$ we are interested in
\begin{equation}
    \De := \sup \left\{\delta \colon\quad
    \begin{matrix}
        \text{every measurable subset of $[0,1)$ with measure $\e$} \\
        \text{contains a symmetric subset with measure $\delta$}
    \end{matrix}
        \right\}.
\label{De.definition}
\end{equation}
It is not immediately obvious, although it turns out to be true, that $\De>0$.

We have dubbed this sort of question ``continuous Ramsey theory'', and we direct the
reader to later in this section for problems with a similar flavor; some of these have appeared in
the literature and some are given here for the first time.

We determine a lower bound for $\De$ using tools and methods from harmonic analysis, some of which were spurred by ideas from nonstandard analysis and the theory of wavelets. We also construct sets without large
symmetric subsets using results from probabilistic
number theory. These two lines of attack complement each other, and our bounds on $\De$ yield new
results in additive number theory as well.

The following theorem, proved in Sections \ref{easy.section} and \ref{deriving.upper.bounds.sec},
states some fundamental properties of the function $\De$.

\begin{thm}
\label{Delta.properties.thm}
The function $\De$ is continuous and, in fact, satisfies the Lipschitz condition
$$|\Delta(x)-\Delta(y)| \le 2|x-y|$$
for all $x,y\in(0,1]$. Furthermore, the function $\frac{\De}{\e^2}$ is increasing on $(0,1]$, and hence
$\lim_{\e\to0^+} \frac{\De}{\e^2}$ exists.
\end{thm}

We turn now to stating our quantitative bounds for $\De$. The lower bound $\De \geq \frac12 \e^2$, which we call
the trivial lower bound on $\De$ (see Lemma~\ref{Delta.Trivial.Lower.Bound.lem} below), is not so far from the
best we can derive. In fact, the bulk of this paper is devoted to improving the constant in this lower bound
from $\frac12$ to $\Dconstant$. Moreover, we are able to establish a complementary upper bound for $\De$ using
results on an analogous problem in combinatorial number theory.

\begin{figure}
    \begin{center}
    \begin{picture}(384,240)
    \put(12,230){$\De/\e^2$}
    \put(370,12){$\e$}
    \put(0,0){\includegraphics{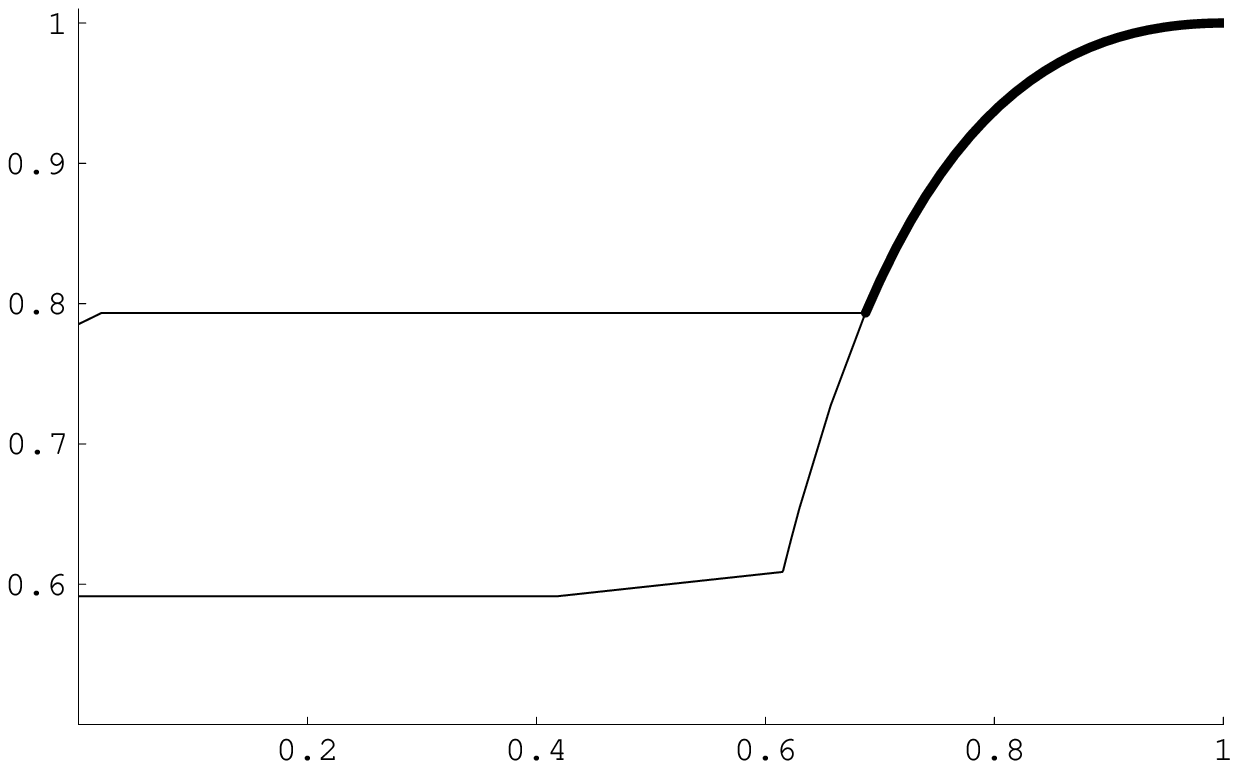}}
    \end{picture}
    \end{center}
    \caption{Upper and Lower Bounds for $\De/\e^2$\label{overallpic}}
\end{figure}

Figure~\ref{overallpic} shows the precise upper and lower bounds we obtain for
$\De/\e^2$ as functions of $\e$, which we present as Theorem~\ref{Delta.Summary.thm}.

\begin{thm}
\label{Delta.Summary.thm} We have:
\renewcommand{\theenumi}{\roman{enumi}}
    \begin{enumerate}
        \item $\Delta(\e)=2\e-1$ for $\frac{11}{16}\le\e\le1$, and $\De \ge2\e-1$ for
                 $\frac12\le\e\le\frac{11}{16}$;
        \item $\De \ge \Dconstant \e^2$ for all $0<\e\le1$;
        \item $\De \ge 0.5546\e^2 + 0.088079\e^3$ for all $0<\e \le1$;
        \item $\De \le \frac{96}{121}\e^2 < 0.79339 \e^2$ for
$0<\e\le\frac{11}{16}$;
        \item $\De \le \frac{\pi\e^2}{(1+\sqrt{1-\e})^2} = \frac\pi4\e^2 +
O(\e^3)$ for all $0<\e\le1$.
    \end{enumerate}
\end{thm}

Note that $\frac\pi4<0.7854$. The upper bound in part (v) of the theorem is superior to
the one in part (iv) in the range $0<\e<\frac{11}{96}(8\sqrt{6\pi}-11\pi) \doteq
0.0201$. The five parts of the theorem are proved separately in
Proposition~\ref{Line.2e-1.prop}, Proposition~\ref{the.full.De.bound.prop},
Proposition~\ref{Delta.one.half.thm}, Corollary~\ref{part.iii.cor}, and
Proposition~\ref{part.iv.prop}, respectively.

Perhaps surprisingly, the upper bounds given in Theorem~\ref{Delta.Summary.thm}(iv)--(v) are
derived from number-theoretic considerations. A set $S$ of integers is called a $B^\ast[g]$ set if
for any given $m$ there are at most $g$ ordered pairs $(s_1,s_2)\in S \times S$ with $s_1+s_2=m$.
We shall use constructions of large $B^\ast[g]$ sets to derive upper bounds on $\De$ in
Section~\ref{deriving.upper.bounds.sec}. Conversely, our bounds on $\De$ improve the best known
upper bounds on the size of $B^\ast[g]$ sets for large~$g$, as we show in \cite{Martin.O'Bryant.b}.
See the article \cite{2004.O'Bryant} of the second author for a survey of $B^\ast[g]$ sets.

We note that the difficulty of determining $\De$ is in stark contrast to the analogous problem where we consider
subsets of the circle $\T:=\R/\Z$ instead of subsets of the interval $[0,1]$. In this analogous setting, we
completely determine the corresponding function $\Delta_{\T}(\e)$; in fact, we show
(Corollary~\ref{circle.answer.cor}) that $\Delta_{\T}(\e)=\e^2$ for all $0<\e\le1$. As it turns out, the methods
that allow the proof of the upper bound $\Delta_{\T}(\e)\le\e^2$, namely constructions of large $B^\ast[g]$ sets
in $\Z/N\Z$, are also helpful to us in constructing the large $B^\ast[g]$ sets themselves.

Schinzel and Schmidt \cite{2002.Schinzel.Schmidt} consider the problem of bounding
    $$B:=\sup_f \frac{\| f\ast f\|_1}{\|f \ast f\|_\infty},$$
where the supremum is taken over all nonnegative functions supported on $[0,1]$; they showed that
$4/\pi \le B < 1.7373$. The proof of Theorem~\ref{Delta.Summary.thm}(ii) improves the value
$1.7373$ to $1.691$.

We remark briefly on the phrase ``continuous Ramsey theory''. A ``coloring theorem'' has the form:
\begin{quote}
Given some fundamental set $R$ colored with a finite number of colors, there exists a highly
structured monochromatic subset, provided that $R$ is sufficiently large.
\end{quote}
The prototypical example is Ramsey's Theorem itself: however one colors the edges of the complete
graph $K_n$ with $r$ colors, there is a monochromatic complete subgraph on $t$ vertices, provided
that $n$ is sufficiently large in terms of $r$ and $t$. Another example is van der Waerden's
Theorem: however one colors the integers $\{1,2,\dots,n\}$ with $r$ colors, there is a
monochromatic arithmetic progression with $t$ terms, provided that $n$ is sufficiently large in
terms of $r$ and $t$.

In many cases, the coloring aspect of a Ramsey-type theorem is a ruse and one may prove
a stronger statement with the form:
\begin{quote}
Given some fundamental set $R$, any large subset of $R$ contains a highly structured subset,
provided that $R$ itself is sufficiently large.
\end{quote}
Such a result is called a ``density theorem.'' For example, van der Waerden's Theorem is a special
case of the density theorem of Szemer\'{e}di: Every subset of $\{1,2,\dots,n\}$ with cardinality at
least $\delta n$ contains a $t$-term arithmetic progression, provided that $n$ is sufficiently
large in terms of $\delta$ and $t$.

Ramsey theory is the study of such theorems on different types of structures. By
``continuous Ramsey theory'' we refer to Ramsey-type problems on continuous measure
spaces. In particular, this paper is concerned with a density-Ramsey problem on the
structure $[0,1)\subseteq\R$ with Lebesgue measure. The type of substructure we focus on
is a symmetric subset.

Other appearances of continuous Ramsey theory in the literature are in the work of
\cite{1958.Swierczkowski} (see also \cite[problem C17]{1994.Guy}),
\cite{2000.Banakh.Verbitsky.Vorobets}, \cite{2002.Schinzel.Schmidt}, and \cite{Chung.Erdos.Graham}.
In all cases, there is an analogous discrete Ramsey theory problem. However, see
\cite{Chung.Erdos.Graham} for an interesting example where the quantities involved in the discrete
setting do not tend in the limit to the analogous quantity in the continuous setting.

\section{Lower Bounds for $\De$}

We give easy lower bounds for $\De$ and prove that $\De$ is continuous in Section~\ref{easy.section}.
Section~\ref{introduction.section} below makes explicit the connection between $\De$ and harmonic
analysis. Section~\ref{basic.argument.section} gives a simple, but quite good, lower bound on
$\De$. In Section~\ref{main.bound.section}, we give a more general form of the argument of
Section~\ref{basic.argument.section}. Using an analytic inequality established in
Section~\ref{useful.inequalities.section}, we investigate in Section~\ref{full.bound.section} the
connection between $\ffi$ and the Fourier coefficients of $f$, which, when combined with the
results of Section~\ref{basic.argument.section}, allows us to show that $\De \geq \Dconstant \e^2$.
The bound given in Section~\ref{basic.argument.section} and improved in
Section~\ref{full.bound.section} depends on a kernel function with certain properties; in
Section~\ref{kernel.problem.section} we discuss how we chose our kernel. In
Section~\ref{delta.half.section}, we use a different approach to derive a lower bound on $\De$
which is superior for $\frac38<\e<\frac58$.

\subsection{Easy Bounds for $\De$}
\label{easy.section}

We now turn our attention to the investigation of the function $\De$
defined in Eq.~(\ref{Deltadef}). In this section we establish several simple
lemmas describing basic properties of $\Delta$.

Let $\lambda$ denote Lebesgue measure on $\R$. We find the following equivalent definition of $\De$ easier to work with than the definition given in Eq.~\eqref{De.definition}: if we define
    \begin{equation}
    D(A) := \sup\{ \lambda(C)\colon\quad C\subseteq A,\, \text{$C$ is symmetric}\},
    \label{Ddef}
    \end{equation}
then
    \begin{equation}
    \De := \inf\{ D(A)\colon\quad A\subseteq [0,1),\, \lambda(A)=\e\}.
    \label{Deltadef}
    \end{equation}

\begin{lem}
\label{Delta.Trivial.Bounds.lem} $\De \ge 2\e-1$ for all $1/2\le\e\le1$.
\end{lem}

\begin{proof}
For $A\subseteq[0,1)$, the centrally symmetric set $A\cap(1-A)$ has measure equal to
    $$
    \lambda(A)+\lambda(1-A)-\lambda(A\cup(1-A))
        = 2\lambda(A)-\lambda(A\cup(1-A)) \geq 2\lambda(A)-1.
    $$
Therefore $D(A)\ge2\lambda(A)-1$ from the definition (\ref{Ddef}) of the function $D$.
Taking the infimum over all subsets $A$ of $[0,1)$ with measure $\e$, this becomes
$\De\ge2\e-1$ as claimed.
\end{proof}

While this bound may seem obvious, it is in many situations the state of the art. As we
show in Proposition~\ref{Line.2e-1.prop} below, $\De$ actually equals $2\e-1$ for
$\tfrac{11}{16}\le\e\le1$; and $\De \geq 2\e-1$ is the best lower bound of which we are
aware in the range $\TrivialLowerBoundconstant \leq \e <\tfrac{11}{16} = 0.6875$.

One is tempted to try to sharpen the bound $\De \geq 2\e-1$ by considering the symmetric
subsets with center $1/3$, $1/2$, or $2/3$, for example, instead of merely $1/2$.
Unfortunately, it can be shown that given any $\e\ge\frac12$ and any finite set
$\{c_1,c_2,\dots,c_n\}$, one can construct a sequence $S_k$ of sets, each with measure
$\e$, that satisfies
    $$
    \lim_{k\to\infty}\left(\max_{1\leq i \leq n}
            \{\lambda\left(S_k\cap(2c_i-S_k)\right)\}\right)
                        = 2\e-1.
    $$
Thus, no improvement is possible with this sort of argument.

\begin{lem}[Trivial Lower Bound]
$\Delta(\e) \ge \frac12\e^2$ for all $0\le\e\le1$.
\label{Delta.Trivial.Lower.Bound.lem}
\end{lem}

\begin{proof}
Given a subset $A$ of $[0,1)$ of measure $\e$, let $A(x)$ denote the indicator function of $A$, so
that the integral of $A(x)$ over the interval $[0,1)$ equals $\e$. If we define $f(c) := \int_{0}^1
A(x) A(2c-x) \, dx$, then $f(c)$ is the measure of the largest symmetric subset of $A$ with center
$c$, and we seek to maximize $f(c)$. But $f$ is clearly supported on $[0,1)$, and so
    \begin{multline*}
    D(A)
    = \max_{0\le c \le 1} f(c) \ge \int_0^1 f(c)\,dc
    = \int_0^1 \int_0^1 A(x) A(2c-x) \, dc \, dx \\
    = \int_0^1 A(x) \bigg(  \int_{-x}^{2-x} A(w)\, \frac{dw}2 \bigg) \, dx
    = \frac 12 \int_0^1 A(x)\,dx\; \int_0^1 A(w)\,dw
    = \tfrac12 \e^2,
    \end{multline*}
since $A(w)$ is supported on $[0,1] \subseteq[-x,2-x]$.
Since $A$ was an arbitrary subset of $[0,1)$ of measure $\e$, we have shown that
\mbox{$\Delta(\e)\ge\frac12\e^2$}.
\end{proof}

It is obvious from the definition of $\Delta$ that $\De$ is an increasing function; the next lemma shows that $\frac{\De}\e$ is also an increasing function.
Later in this paper (see Proposition~\ref{De.over.e2.increasing.prop}), we will show that in fact
even $\frac{\De}{\e^2}$ is an increasing function.

\begin{lem}
\label{Delta.Line.To.Origin.lem} $\Delta(\e) \leq \frac{\Delta(x)}{x}\e$ for all $0\leq
\e \leq x\le1$. In particular, $\De \le \e$.
\end{lem}

\begin{proof}
If $tA := \{ta\colon a\in A\}$ is a scaled copy of a set $A$, then clearly $D(tA)=tD(A)$. Applying
this with any set $A\subseteq[0,1)$ of measure $x$ and with $t=\frac\e{x}\le1$, we see that
$\frac\e{x}A$ is a subset of $[0,1)$ with measure $\e$, and so by the definition of $\Delta$ we
have $\Delta(\e) \le D(\frac\e{x}A) = \frac\e{x}D(A)$. Taking the infimum over all sets
$A\subseteq[0,1)$ of measure $x$, we conclude that $\Delta(\e) \le \frac\e{x}\Delta(x)$. The second
assertion of the lemma follows from the first assertion with the trivial value $\Delta(1)=1$.
\end{proof}

Let $$S\oplus T := (S\setminus T)\cup(T\setminus S)$$ denote the symmetric difference of $S$ and
$T$. (While this operation is more commonly denoted with a triangle rather than with a $\oplus$, we
would rather avoid any potential confusion with the function $\Delta$ featured prominently in this
paper.)

\begin{lem}
If $S$ and $T$ are two sets of real numbers, then $|D(S)-D(T)| \le 2\lambda(S\oplus T)$.
\label{Diamond.lem}
\end{lem}

\begin{proof}
Let $E$ be any symmetric subset of $S$, and let $c$ be the center of $E$, so that $E=2c-E$. Define
$F=E\cap T\cap(2c-T)$, which is a symmetric subset of $T$ with center $c$. We can write
$\lambda(F)$ using the inclusion-exclusion formula
    \begin{multline*}
    \lambda(F) = \lambda(E) + \lambda(T) + \lambda(2c-T) \\ - \lambda(E\cup T) -
    \lambda(E\cup(2c-T)) - \lambda(T\cup(2c-T)) + \lambda(E\cup T\cup(2c-T)).
    \end{multline*}
Rearranging terms, and noting that $T\cup(2c-T) \subseteq E\cup T\cup(2c-T)$, we see
that
    \begin{equation*}
    \lambda(E) - \lambda(F) \le -\lambda(T)-\lambda(2c-T)+\lambda(E\cup T)
            +\lambda(E\cup (2c-T)).
    \end{equation*}
Because reflecting a set in the point $c$ does not change its measure, this is
the same as
    \begin{equation*}
    \begin{split}
    \lambda(E) - \lambda(F) &\le -\lambda(T)-\lambda(T)+\lambda(E\cup T)
            +\lambda(E\cup (2c-T)) \\
    &= -\lambda(T)-\lambda(T)+\lambda(E\cup T)+\lambda((2c-E)\cup T) \\
    &= 2\big( \lambda(E\cup T) - \lambda(T) \big) \\
    &\le 2\big( \lambda(S\cup T) - \lambda(T) \big)
    = 2\lambda(S\setminus T) \le 2\lambda( S\oplus T).
    \end{split}
    \end{equation*}
Therefore, since $F$ is a symmetric subset of $T$,
\begin{equation*}
\lambda(E) \le \lambda(F) + 2\lambda( S\oplus T) \le D(T) + 2\lambda( S\oplus T).
\end{equation*}
Taking the supremum over all symmetric subsets $E$ of $S$, we conclude that $D(S) \le D(T) +
2\lambda( S\oplus T)$. If we now exchange the roles of $S$ and $T$, we see that the proof is
complete.
\end{proof}

\begin{lem}
The function $\Delta$ satisfies the Lipschitz condition
$|\Delta(x)-\Delta(y)|\leq 2|x-y|$ for all $x$ and $y$ in $[0,1]$. In
particular, $\Delta$ is continuous. \label{Lipschitz.Condition.lem}
\end{lem}

\begin{proof}
Without loss of generality assume $y<x$. In light of the monotonicity $\Delta(y)\le \Delta(x)$, it
suffices to show that $\Delta(y) \ge \Delta(x) - 2(x-y)$. Let $S\subseteq[0,1)$ have measure $y$.
Choose any set $R\subseteq[0,1)\setminus S$ with measure $x-y$, and set $T=S\cup R$. Then $S\oplus
T=R$, and so by Lemma~\ref{Diamond.lem}, $D(T)-D(S) \le 2\lambda(R) = 2(x-y)$. Therefore $D(S) \ge
D(T) - 2(x-y) \ge \Delta(x) - 2(x-y)$ by the definition of $\Delta$. Taking the infimum over all
sets $S\subseteq[0,1)$ of measure $y$ yields $\Delta(y) \ge \Delta(x) - 2(x-y)$ as desired.
\end{proof}

\subsection{Notation}\label{introduction.section}

There are many ways to define the basic objects of Fourier analysis; we follow
\cite{1984.Folland}. Unless specifically noted otherwise, all integrals are over the
circle group $\T:=\R/\Z$; for example, $L^1$ denotes the class of functions $f$ for
which $\int_\T |f(x)|\,dx$ is finite.
For each integer $j$, we define $\hat{f}(j):=\int
f(x) e^{-2\pi ijx}\,dx$, so that for any function $f\in L^1$, we have
$f(x)=\sum_{j=-\infty}^\infty \hat{f}(j)e^{2\pi ij x}$ almost everywhere. We define the
convolution $f*g(c):=\int f(x)g(c-x) \,dx$, and we note that $\widehat{f\ast
g}(j)=\hat{f}(j)\hat g(j)$ for every integer $j$; in particular,
$\widehat{\ff}(j)=\hat{f}(j)^2$.

We define the usual $L^p$ norms $$\textstyle\|f\|_p := \big( \int |f(x)|^p\,
dx\big)^{1/p}$$ and
    $$
   \|f\|_\infty    := \lim_{p\to\infty} \|f\|_p
                    = \sup\big\{y\colon \lambda(\{x\colon |f(x)|>y\})>0\big\}.
    $$
With these definitions, H\"{o}lder's Inequality is valid: if $p$ and $q$ are conjugate
ex\-po\-nents---that is, $\frac1p+\frac1q=1$---then $\|fg\|_1 \leq \|f\|_p \|g\|_q$. We
also note that $\|f*g\|_1 = \|f\|_1\|g\|_1$ when $f$ and $g$ are nonnegative functions; in particular,
$\|\ff\|_1=\|f\|_1^2=\hat{f}(0)^2$. We shall also employ the $\ell^p$ norms for
bi-infinite sequences: if $a=\{a_j\}_{j\in\Z}$, then $\|a\|_p := \big( \sum_{j\in\Z}
|a_j|^p \big)^{1/p}$ and $\|a\|_\infty := \lim_{p\to\infty} \|a\|_p = \sup_{j\in\Z}
|a_j|$. Although we use the same notation for the $L^p$ and $\ell^p$ norms, no confusion
should arise, as the object inside the norm symbol will either be a function on $\T$ or
its sequence of Fourier coefficients. With this notation, we recall
Parseval's identity $$\int f(x)g(x)\,dx = \sum \hat f(j)\hat g(-j)$$ (assuming the
integral and sum both converge); in particular, if $f=g$ is real-valued (so that $\hat
f(-j)$ is the conjugate of $\hat f(j)$ for all $j$), this becomes $\|f\|_2 =
\|\hat{f}\|_2$. The Hausdorff-Young inequality, $\|\hat{f}\|_q \leq \| f \|_p$ whenever
$p$ and $q$ are conjugate exponents with $1\le p\le 2\le q\le\infty$, can be thought of
as a generalization of this latter version of Parseval's identity. We also require the
definition
    \begin{equation}
    \lnorm{m}{p}{a} = \left(\sum_{|j|\geq m} |a(j)|^p \right)^{1/p}
    \label{lnorm.def}
    \end{equation}
for any sequence $a=\{a_j\}_{j\in\Z}$, so that $\lnorm0pa = \|a\|_p$, for example.

We note that for any fixed sequence $a=\{a_j\}_{j\in\Z}$, the $\ell^p$-norm $\|a\|_p$ is a
decreasing function of $p$. To see this, suppose that $1\le p\le q<\infty$ and
$a\in\ell^p$. Then $|a_j| \le \|a\|_p$ for all $j\in\Z$, whence $|a_j|^{q-p} \le
\|a\|_p^{q-p}$ (since $q-p\ge0$) and so $|a_j|^q \le \|a\|_p^{q-p} |a_j|^p$. Summing
both sides over all $j\in\Z$ yields $\|a\|_q^q \le \|a\|_p^{q-p} \|a\|_p^p = \|a\|_p^q$,
and taking $q$th roots gives the desired inequality $\|a\|_q \le \|a\|_p$.

Finally, we define a ``pdf'', short for ``probability density function'', to be a
nonnegative function in $L^2$ whose $L^1$-norm (which is necessarily finite, since $\T$
is a finite measure space) equals~1. Also, we single out a special type of pdf called an
``nif'', short for ``normalized indicator function'', which is a pdf that only takes one
nonzero value, that value necessarily being the reciprocal of the measure of the support
of the function. (We exclude the possibility that an nif takes the value 0 almost
everywhere.) Specifically, we define for each $E\subseteq\T$ the nif
    $$
    f_E(x):=
        \begin{cases}
            \lambda(E)^{-1} & x\in E, \\
            0 & x\not\in E.
        \end{cases}
    $$
Note that if $f$ is a pdf, then $1=\hat{f}(0)=\hat{f}(0)^2=\|f\|_1^2=\|\ff\|_1$.

We are now ready to reformulate the function $\De$ in terms of this notation.

\begin{lem}
\label{Fourier.Connection.lem} We have
    $$ \textstyle
    \frac12 \e^2 \inf_g \|g*g\|_\infty \le \frac12 \e^2 \inf_f \ffi = \De,
    $$
the first infimum being taken over all pdfs $g$ that are supported on
$[-\tfrac14,\tfrac14]$, and the second infimum being taken over all nifs $f$ whose
support is a subset of $[-\frac14,\frac14]$ of measure $\frac\e2$.
\end{lem}

\begin{proof}
The inequality is trivial, since every nif is a pdf; it remains to prove the equality.

For each measurable $A\subseteq[0,1)$, define $E_A:=\{\frac12 (a-\frac12) \colon a\in
A\}\subseteq[-\frac14,\frac14]$. The sets $A$ and $E_A$ differ only by translation and
scaling, so that $\lambda(A)=2\lambda(E_A)$ and $D(A)=2D(E_A)$. Thus
    \begin{align*}
    \De     &:= \inf\{ D(A) \colon A\subseteq[0,1),\, \lambda(A)=\e\} \\
            &=  \e^2 \inf\left\{ \frac{D(A)}{\lambda(A)^2}
                            \colon A\subseteq[0,1),\, \lambda(A)=\e\right\} \\
            &=  \e^2 \inf\left\{ \frac{2D(E_A)}{(2\lambda(E_A))^2}
                            \colon A\subseteq[0,1),\, \lambda(A)=\e\right\} \\
            &=  \frac12 \e^2 \inf\left\{ \frac{D(E)}{\lambda(E)^2}
                            \colon E\subseteq[-\tfrac14,\tfrac14],\,
                            \lambda(E)=\tfrac{\e}{2}              \right\}.
    \end{align*}
For each $E\subseteq[-\tfrac14,\tfrac14]$ with $\lambda(E)=\frac{\e}{2}$, the function
$f_E(x)$ is an nif supported on a subset of $[-\tfrac14,\tfrac14]$ with measure
$\frac{\e}{2}$, and it is clear that every such nif arises from some set $E$. Thus, it
remains only to show that $\frac{D(E)}{\lambda(E)^2}=\| f_E \ast f_E \|_\infty$, i.e.,
that $D(E)=\lambda(E)^2 \| f_E \ast f_E \|_\infty$.

Fix $E\subseteq[-\tfrac14,\tfrac14]$, and let $E(x)$ be the indicator function of $E$.
Note that $f_E(x) = \lambda(E)^{-1}E(x)$. The maximal symmetric subset of $E$ with
center $c$ is $E\cap(2c-E)$, and this has measure $\int E(x) E(2c-x)\,dx$. Thus
    \begin{align*}
    D(E)    &:=\sup\{ \lambda(C)\colon\quad C\subseteq E,\, \text{$C$ is symmetric}\}\\
            &= \sup_c \left( \int E(x)E(2c-x)\,dx \right)\\
            &= \sup_c \left(
                    \int\lambda(E)f_E(x)\lambda(E)f_E(2c-x)\,dx \right)\\
            &= \lambda(E)^2 \sup_c \left( \int f_E(x)f_E(2c-x)\,dx \right)\\
            &= \lambda(E)^2 \sup_c f_E \ast f_E (2c) \\
            &= \lambda(E)^2 \left\| f_E \ast f_E \right\|_\infty,
    \end{align*}
as desired.
\end{proof}

The convolution in Lemma~\ref{Fourier.Connection.lem} may be taken over $\R$ or over
$\T$, the two settings being equivalent since $\ff$ is supported on an interval of
length 1. In fact, the reason we scale $f$ to be supported on an interval of length
$1/2$ is so that we may replace convolution over $\R$, which is the natural place to
study $\De$, with convolution over $\T$, which is the natural place to do harmonic
analysis.

\subsection{The Basic Argument}\label{basic.argument.section}

We begin the process of improving upon the trivial lower bound for $\De$ by stating a
simple version of our method that illustrates the ideas and techniques involved.

\begin{prop}
\label{Cf.First.Result.prop} Let $K$ be any continuous function on $\T$ satisfying
$K(x)\geq1$ when $x\in[-\tfrac14,\tfrac14]$, and let $f$ be a pdf supported on
$[-\tfrac14,\tfrac14]$. Then $$\ffi \geq \|\ff\|_2^2 \geq \|\hat{K}\|_{4/3}^{-4}.$$
\end{prop}

\begin{proof}
We have
$$
1=\int f(x)\,dx \leq \int f(x)K(x)\,dx = \sum_j \hat{f}(j)\hat{K}(-j)
$$
by Parseval's identity. H\"{o}lder's Inequality now gives $1\leq \|\hat{f}\|_4
\|\hat{K}\|_{4/3}$, which we restate as the inequality $\|\hat{K}\|_{4/3}^{-4} \leq
\|\hat{f}\|_4^4$.

Now $\|\hat{f}\|_4^4 = \sum_j |\hat{f}(j)|^4=\sum_j |\widehat{\ff}(j)|^2 = \|\ff\|_2^2$
by another application of Parseval's identity. Since $(\ff)^2\leq \ffi(\ff)$,
integration yields $\|\ff\|_2^2 \leq \ffi\|\ff\|_1 = \ffi$. Combining the last three
sentences, we see that
    $  \|\hat{K}\|_{4/3}^{-4}
                \leq \|\hat{f}\|_4^4
                = \| \ff \|_2^2
                \leq \ffi$
as claimed.
\end{proof}

This reasonably simple theorem already allows us to give a nontrivial lower bound
for~$\De$. The step function
    $$
    K_1(x) :=
        \begin{cases}
            1   &   0\le|x|\le \frac14,\\
            1 - \frac{2\pi^4}{\pi^4 + 24{\zeta(\frac{4}{3})}^3
                \left( 5 + 2^{4/3} - 2^{8/3} \right)}
                     &   \frac14<|x|\le \frac12,
        \end{cases}
    $$
has $\|\hat{K_1}\|_{4/3}^{-4} = 1 + \frac{\pi^4} {8\,{\left( 2^{4/3} -1\right) }^3\,
{\zeta(\frac{4}{3})}^3} > 1.074 $ (the elaborate constant used in the definition of
$K_1$ was chosen to minimize $\|\hat{K_1}\|_{4/3}$). A careful reader may complain that
$K_1$ is not continuous. The continuity condition is not essential, however, as we may
approximate $K_1$ by a continuous function $L$ with $\|L\|_{4/3}$ arbitrarily close to
$\|K_1\|_{4/3}$.

Green \cite{2001.Green} used a discretization of the kernel function
    $$
    K_2(x) :=
        \begin{cases}
            1   &   0\le|x|\le \frac14,\\
            1 - \alpha  +  \alpha \left( 40 (2x - 1)^4 - \frac 32 \right)
                     &   \frac14<|x|\le \frac12,
        \end{cases}
    $$
with a suitably chosen $\alpha$ to get $\|\hat{K_2}\|_{4/3}^{-4}> \frac87>1.142$. We get
a slightly larger value of $\|\hat{K}\|_{4/3}^{-4}$ in the following corollary with a
much more complicated kernel. See Section~\ref{kernel.problem.section} for a discussion
of how we came to find our kernel.

\begin{cor}\label{First.cor}
If $f$ is a pdf supported on $[-\frac14,\frac14]$, then
    $$\|\ff\|_2^2 \geq \fstarftwonormconstant.$$
Consequently, $\De \geq 0.574575 \e^2$ for all $0\le\e\le1$.
\end{cor}

\begin{proof}
Set
    \begin{equation}
    K_3(x):=
    \begin{cases}
        1           &   0\leq |x| \leq \tfrac14, \\
        0.6644+0.3356\left(\tfrac{2}{\pi}
        \tan^{-1}\left(\tfrac{1-2x}{\sqrt{4x-1}}\right)\right)^{1.2015}
                    &   \tfrac14 \le |x| \le \tfrac12.
    \end{cases}
    \label{easy.K.def}
    \end{equation}
$K_3(x)$ is pictured in Figure~\ref{GoodK}.
    \begin{figure}
    \begin{center}
    \begin{picture}(384,128)
    \put(176,114){$K_3(x)$}
    \put(356,0){$x$}
    \put(0,-10){\includegraphics{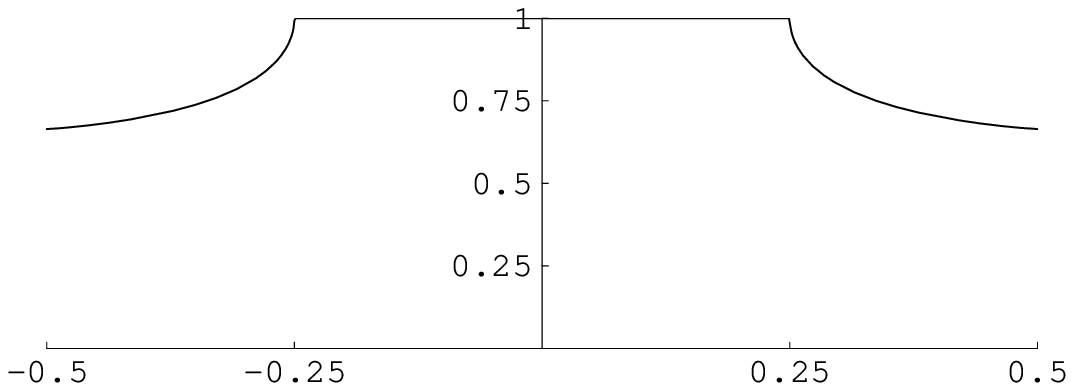}}
    \end{picture}
    \end{center}
    \caption{The function $K_3(x)$\label{GoodK}}
    \end{figure}

We do not know how to rigorously bound $\|\hat{K_3}\|_{4/3}$, but we can rigorously bound $\|\hat{K_4}\|_{4/3}$
where $K_4$ is a piecewise linear function `close' to $K_3$. Specifically, let $K_4(x)$ be the even piecewise
linear function with corners at
    $$(0,1),\left(\frac14,1\right), \left(\frac14+\frac{t}{4\times 10^4},
    K_3\Big(\frac14+\frac{t}{4\times10^4}\Big)\right) \quad (t=0,1,\dots,10^4).$$
We calculate (using Proposition~\ref{Piecewise.Linear.prop} below) that
$\|\hat{K_4}\|_{4/3} < 0.9658413$. Therefore, by Proposition~\ref{Cf.First.Result.prop}
we have
    $$\|\ff\|_2^2 \geq (0.9658413)^{-4} > \fstarftwonormconstant.$$
Using Lemma~\ref{Fourier.Connection.lem}, we now have $\De >
\frac12\e^2(\fstarftwonormconstant) > 0.574575 \e^2$.
\end{proof}

The constants in the definition \eqref{easy.K.def} of $K_3(x)$ were numerically
optimized to minimize $\|\hat{K_4}\|_{4/3}$ and otherwise have no special significance.
The definition of $K_3(x)$ is certainly not obvious, and there are much simpler kernels that
do give nontrivial bounds. In Section~\ref{kernel.problem.section} below, we indicate
the experiments that led to our choice.

We note that the function
    $$
    b(x):=
        \left\{%
        \begin{array}{ll}
            \frac{4/\pi}{\sqrt{1-16x^2}} & -1/4<x<1/4, \\
            0 & \hbox{otherwise}, \\
        \end{array}%
        \right.
    $$
has $\int b=1$ and $\|b\ast b\|_2^2 <1.14939$. Although $b$ is not a pdf (it is not in
$L^2$), it provides strong evidence that the bound on $\|\ff\|_2^2$ given in
Corollary~\ref{First.cor} is not far from best possible.


This bound on $\|\ff\|_2^2$ may be nearly correct, but the resulting bound on $\ffi$ is
not: we prove below that $\ffi \geq \twotimesDconstant$, and believe that $\ffi \geq
\pi/2$. We have tried to improve the argument given in
Proposition~\ref{Cf.First.Result.prop} in the following four ways:
\begin{enumerate}
\item Instead of considering the sum $\sum_{j} \hat{f}(j)\hat{K}(-j)$ as a whole, we
    separate the central terms from the tails and establish inequalities that depend upon the two in distinct ways.
    This generalized form of the above argument is expounded in the next section. The
    success of this generalization relies on certain inequalities restricting the possible
    values of these central coefficients; establishing these restrictions is the goal of
    Sections~\ref{useful.inequalities.section} and~\ref{full.bound.section}. The
    final lower bound derived from these methods is given in
    Section~\ref{full.bound.section}.
\item We have searched for more advantageous kernel functions $K(x)$ for which we can
    compute $\|\hat{K}\|_{4/3}$ in an accurate way. A detailed discussion of our search for
    the best kernel functions is in Section~\ref{kernel.problem.section}.
\item The application of Parseval's identity can be replaced with the Hausdorff-Young
    inequality, which leads to the conclusion $\|\ff\|_\infty \geq \|\hat{K}\|_p^{-q}$,
    where $p\leq \frac{4}{3}$ and $q\geq4$ are conjugate exponents. Numerically, the values
    $(p,q)=(\frac{4}{3},4)$ appear to be optimal. However, Beckner's sharpening
    \cite{Beckner} of the Hausdorff-Young inequality leads to the stronger conclusion
    $\|\ff\|_\infty \geq C(q)\|\hat{K}\|_p^{-q}$ where $C(q) = \frac q2(1-\frac2q)^{q/2-1} =
    \frac q{2e}+O(1)$. We have not experimented to see whether a larger lower bound can be
    obtained from this stronger inequality by taking $q>4$.
\item Notice that we used the inequality $\| g \|_2^2 \leq \|g\|_\infty \|g\|_1$ with
    the function $g=\ff$. This inequality is sharp exactly when the function $g$ takes only
    one nonzero value (i.e., when $g$ is a nif), but the convolution $\ff$ never
    behaves that way. Perhaps for these autoconvolutions, an analogous inequality with a
    stronger constant than 1 could be established. Unfortunately, we have not been able to
    realize any success with this idea, although we believe
    Conjecture~\ref{infinity.two.cnj} below. If true, the conjecture implies
    the bound $\De \geq 0.651 \e^2$.
\end{enumerate}
    \begin{cnj}\label{infinity.two.cnj}
    If $f$ is a pdf supported on $[-\frac14,\frac14]$, then
        $$\frac{\ffi }{\|\ff\|_2^2} \geq \frac{\pi}{\log 16},$$
    with equality only if either $f(x)$ or $f(-x)$ equals $\sqrt{\frac2{4x+1}}$ on the
    interval $|x|\le\frac14$.
    \end{cnj}

We remark that Proposition~\ref{Cf.First.Result.prop} can be extended from a twofold
convolution in one dimension to an $h$-fold convolution in $d$ dimensions.

\begin{prop}
Let $K$ be any continuous function on $\T^d$ satisfying $K(\bar{x})\geq1$ when
$\bar{x}\in[-\tfrac1{2h},\tfrac1{2h}]^d$, and let $f$ be a pdf supported on
$[-\tfrac1{2h},\tfrac1{2h}]^d$. Then
    $$\| f^{\ast h}\|_\infty \geq \|f^{\ast h}\|_2^2 \geq \|\hat{K}\|_{2h/(2h-1)}^{-2h}.$$
Every subset of $[0,1]^d$ with measure $\e$ contains a symmetric subset with measure
$(0.574575)^d \e^2$.
\end{prop}

\begin{proof}
The proof proceeds as above, with the conjugate exponents $(\frac{2h}{2h-1},2h)$ in
place of $(\frac43,4)$, and the kernel function $K(x_1,x_2,\dots,x_d)=K(x_1)K(x_2)\cdots
K(x_d)$ in place of the kernel function $K(x)$ defined in the proof of
Corollary~\ref{First.cor}. The second assertion of the proposition follows on taking $h=2$.
\end{proof}

\subsection{The Main Bound}\label{main.bound.section}

We now present a more subtle version of Proposition \ref{Cf.First.Result.prop}. Recall
that the notation $\lnorm npa$ was defined in Eq.~(\ref{lnorm.def}). We also use $\Re z$
to denote the real part of the complex number~$z$.

\begin{prop}
\label{Main.Bound.prop} Let $m\ge1$. Suppose that $f$ is a
pdf supported on $[-\tfrac14,\tfrac14]$ and that $K$ is even, continuous, satisfies
$K(x) = 1$ for $-\tfrac14 \le x <\frac14$, and $\lnorm{m}{4/3}{\hat{K}}>0$. Set $M := 1-
\hat{K}(0) - 2\sum_{j=1}^{m-1} \hat{K}(j) \Re\hat{f}(j)$. Then
\begin{equation}\label{qfinite}
  \|\ff\|_2^2 = \sum_{j\in\Z} |\hat{f}(j)|^4 \geq
            1+\bigg(\frac{M}{\lnorm{m}{4/3}{\hat{K}}}\bigg)^4
            +2 \sum_{j=1}^{m-1} |\Re\hat{f}(j)|^4.
\end{equation}
\end{prop}

\begin{proof}
The equality follows from Parseval's formula
    $$\|\ff\|_2^2 = \sum_j |\widehat{\ff}(j)|^2 = \sum_j |\hat{f}(j)|^4.$$
As in the proof of Proposition~\ref{Cf.First.Result.prop}, we have
    $$1         =   \int f(x)K(x)\,dx
                =   \sum_j \hat{f}(j)\hat{K}(-j)
                =   \sum_{|j|<m} \hat{f}(j)\hat{K}(-j) +
                         \sum_{|j|\geq m} \hat{f}(j)\hat{K}(-j).$$
Since $K$ is even, $\hat{K}(-j)=\hat{K}(j)$ is real, and since $f$ is real valued,
$\hat{f}(-j)=\overline{\hat{f}(j)}$. We have
    $$1      =   \hat{K}(0)+2\sum_{j=1}^{m-1} \hat{K}(j) \Re\hat{f}(j)
                    +\sum_{|j|\geq m} \hat{f}(j) \hat{K}(j),$$
which we can also write as $M=\sum_{|j|\geq m} \hat{f}(j)\hat{K}(j)$. Taking absolute values and applying H\"{o}lder's inequality, we have
    \begin{align*}
    |M|     &\leq \sum_{|j|\geq m} |\hat{f}(j)\hat{K}(j)| \\
            &\leq \left(\sum_{|j|\geq m} |\hat{f}(j)|^4\right)^{1/4}
                 \left(\sum_{|j|\geq m} |\hat{K}(j)|^{4/3}\right)^{3/4} \\
            &= \left(\sum_{|j|\geq m} |\hat{f}(j)|^4\right)^{1/4}
                    \lnorm{m}{4/3}{\hat{K}},
    \end{align*}
which we recast in the form
    $$\sum_{|j|\geq m} |\hat{f}(j)|^4 \geq \bigg(\frac{M}{\lnorm{m}{4/3}{\hat{K}}} \bigg)^4.$$
We add $\sum_{|j|<m} |\hat{f}(j)|^4$ to both sides and observe that $\hat{f}(0)=1$ and
$|\hat{f}(j)| \geq |\Re \hat{f}(j)|$ to finish the proof of the inequality.
\end{proof}

With $x_j:=\Re \hat{f}(j)$, the bound of
Proposition~\ref{Main.Bound.prop} becomes
    $$ 1+\bigg(\frac{1-\hat{K}(0)-2\sum_{j=1}^{m-1}\hat{K}(j)x_j}
           {\lnorm{m}{4/3}{\hat{K}}} \bigg)^4 + 2 \sum_{j=1}^{m-1} x_j^4.$$
This is a quartic polynomial in the $x_j$, and consequently it is not difficult to
minimize, giving an absolute lower bound on $\| \ff \|_2^2$. This minimum occurs at
    $$x_j=\frac{(\hat{K}(j))^{1/3}
    \left(1-\hat{K}(0)-
    2\sum_{i=1}^{j-1}\hat{K}(i)x_i\right)}{\lnorm{j}{4/3}{\hat{K}}^{4/3}},$$
where $(\hat{K}(j))^{1/3}$ is the real cube root of $\hat{K}(j)$. A substitution and simplification of the resulting expression then yields
 $$\min_{x_j\in\R} \left\{1+\bigg(\frac{1-\hat{K}(0)-2\sum_{j=1}^{m-1}\hat{K}(j)x_j}
                    {\lnorm{m}{4/3}{\hat{K}}}     \bigg)^4
            +2 \sum_{j=1}^{m-1} x_j^4 \right\}=
            1+\left( \frac{1-\hat{K}(0)}{\lnorm{1}{4/3}{\hat{K}}} \right)^4,$$
which is nothing more than the bound that Proposition~\ref{Main.Bound.prop} gives with
$m=1$. Moreover,
    $$ 1+\left( \frac{1-\hat{K}(0)}{\lnorm{1}{4/3}{\hat{K}}} \right)^4
        = \sup_{0\le \alpha \le 1} \| (\alpha+(1-\alpha)K)^{\wedge} \|_{4/3}^{-4}, $$
(the details of this calculation are given in Section~\ref{kernel.problem.section}) so that
Proposition~\ref{Main.Bound.prop}, by itself, does not give a different bound on $\| \ff
\|_2^2$ than Proposition~\ref{Cf.First.Result.prop}.

However, we shall obtain additional information on $\hat{f}(j)$ in terms of $\ffi$ in
Section~\ref{full.bound.section} below, and this information can be
combined with Proposition~\ref{Main.Bound.prop} to provide a stronger lower bound on
$\ffi$ than that given by Proposition~\ref{Cf.First.Result.prop}.

%

\begin{cor}
    \label{Cf.Main.Bound.cor}
Let $f$ be a pdf supported on $[-\frac14,\frac14]$, and set $x_1:=\Re\hat{f}(1)$. Then
    \begin{equation*}
    \|\ff\|_2^2 \geq \sum_{j\in\Z} |\hat{f}(j)|^4
        \ge 1 +2 x_1^4
        +(1.53890149 - 2.26425375 x_1)^4.
    \end{equation*}
\end{cor}

\begin{proof}
Set
    $$
    K_5(x)=
    \begin{cases}
        1 & |x|\leq\tfrac14, \\
        1 - (1 - (4(\tfrac12 - x))^{1.61707})^{0.546335} & \tfrac14 < |x|\leq \tfrac12.
    \end{cases}
    $$
Denote by $K_6(x)$ the even piecewise linear function with corners at
 $$
 (0,1), \left(\frac14,1\right),\left(\frac14+\frac{t}{4\times10^4},
    K_5\bigg(\frac14+\frac{t}{4\times10^4}\bigg)\right) \quad(t=1,\dots,10^4).
    $$
We find (using Proposition~\ref{Piecewise.Linear.prop}) that $\hat{K_6}(0)\doteq
0.631932628$, $\hat{K_6}(1)\doteq 0.270776892$, and $\lnorm{2}{4/3}{\hat{K_6}}\doteq
0.239175395$. Apply Proposition~\ref{Main.Bound.prop} with $m=2$ to finish the proof.
\end{proof}

\subsection{Some Useful Inequalities}\label{useful.inequalities.section}

Hardy, Littlewood, and P\'{o}lya~\cite{1988.Hardy.Littlewood.Polya} call a function $u(x)$
{\it symmetric decreasing\/} if $u(x)=u(-x)$ and $u(x) \geq u(y)$ for all $0\leq x\leq
y$, and they call
 $$\sdr{f}(x):= \inf\left\{ y \colon
      \lambda\left(\left\{t \colon f(t)\geq y \right\}\right)\le 2|x| \right\}$$
the {\it symmetric decreasing rearrangement\/} of $f$. For example, if $f$ is the
indicator function of a set with measure $\mu$, then $\sdr{f}$ is simply the indicator
function of the interval $(-\frac\mu2,\frac\mu2)$. Another example is any function $f$
defined on an interval $[-a,a]$ and is periodic with period $\frac{2a}n$, where $n$ is a
positive integer, and that is symmetric decreasing on the subinterval $[-\frac an,\frac
an]$; then $\sdr{f}(x) = f(\frac xn)$ for all $x\in[-a,a]$. In particular, on the
interval $[-\frac14,\frac14]$, we have $\sdr{\cos}(2\pi j x)=\cos(2\pi x)$ for any
nonzero integer $j$. We shall need the following result \cite[Theorem 378]{1988.Hardy.Littlewood.Polya}:
$$
\int f(x)u(x)\,dx \leq \int \sdr{f}(x)\sdr{u}(x)\,dx.
\label{1988.Hardy.Littlewood.Polya.ineq}
$$

We say that {\em $\bar{f}$ is more focused than $f$} (and $f$ is less focused than
$\bar{f}$) if for all $z\in[0,\frac12]$ and all $r\in\T$ we have
 $$\int_{r-z}^{r+z} f \leq \int_{-z}^z \bar{f}.$$
For example, $\sdr{f}$ is more focused than $f$. In fact, we introduce this terminology
because it refines the notion of symmetric decreasing rearrangement in a way that is
useful for us. To give another example, if $f$ is a nonnegative function, set $\bar f$
to be $\|f\|_\infty$ times the indicator function of the interval
$[-\frac1{2\|f\|_\infty},\frac1{2\|f\|_\infty}]$; then $\bar f$ is more focused than
$f$.

\begin{lem}
    \label{Focus.lem}
Let $u(x)$ be a symmetric decreasing function, and let
$h,\bar{h}$ be pdfs with $\bar{h}$ more focused than $h$. Then for all $r\in\T$,
    $$ \int h(x-r) u(x) \, dx \leq \int \bar{h}(x) u(x) \, dx.$$
\end{lem}

\begin{proof}
Without loss of generality we may assume that $r=0$, since if $\bar h(x)$ is more
focused than $h(x)$, then it is also more focused than $h(x-r)$. Also, without loss of
generality we may assume that $h,\bar h$ are continuous and strictly positive on $\T$,
since any nonnegative function in $L^1$ can be $L^1$-approximated by such.

Define $H(z) = \int_{-z}^z h(t)\,dt$ and $\bar H(z) = \int_{-z}^z \bar h(t)\,dt$, so
that $H(\frac12) = \bar H(\frac12) = 1$, and note that the more-focused hypothesis
implies that $H(z)\le\bar H(z)$ for all $z\in[0,\frac12]$. Now $h$ is continuous and
strictly positive, which implies that $H$ is differentiable and strictly increasing on
$[0,\frac12]$ since $H'(z)=h(z)+h(-z)$. Therefore $H^{-1}$ exists as a function from
$[0,1]$ to $[0,\frac12]$. Similar comments hold for $\bar H^{-1}$.

Since $H\le\bar H$, we see that $\bar H^{-1}(s) \le H^{-1}(s)$ for all $s\in[0,1]$.
Then, since $H^{-1}(s)$ and $H^{-1}(s)$ are positive and $u$ is decreasing for positive
arguments, we conclude that $u(H^{-1}(s)) \le u(\bar H^{-1}(s))$, and so
    \begin{equation}
    \int_0^1 u(H^{-1}(s))\,ds \le \int_0^1 u(\bar H^{-1}(s))\,ds.
    \label{H.inverse.inequality}
    \end{equation}
On the other hand, making the change of variables $s=H(t)$, we see that
    $$
    \int_0^1 u(H^{-1}(s))\,ds = \int_0^{H^{-1}(1)} u(t) H'(t)\,dt = \int_0^{1/2}
    u(t) (h(t)+h(-t))\,dt = \int_{\T} u(t)h(t)\,dt
    $$
since $u$ is symmetric. Similarly $\int_0^1 u(\bar H^{-1}(s))\,ds = \int_{\T} u(t)\bar
h(t)\,dt$, and so inequality (\ref{H.inverse.inequality}) becomes $\int u(t)h(t)\,dt \le
\int u(t)\bar h(t)\,dt$ as desired.
\end{proof}

\subsection{The Full Bound}\label{full.bound.section}

To use Proposition~\ref{Main.Bound.prop} to bound $\De$, we need to develop a better
understanding of the central Fourier coefficients $\hat{f}(j)$ for small $j$. In
particular, we wish to apply Proposition~\ref{Main.Bound.prop} with $m=2$, i.e., we need
to develop the connections between $\ffi$ and the real part of the Fourier coefficient
$\hat{f}(1)$.

We turn now to bounding $|\hat{f}(j)|$ in terms of $\ffi$. The guiding principle is that
if $\ff$ is very concentrated then $\ffi$ will be large, and if $\ff$ is not very
concentrated then $|\hat{f}(j)|$ will be small. Green \cite[Lemma 26]{2001.Green} proves
the following lemma in a discrete setting, but since we need a continuous version we
include a complete proof.

\begin{lem}
\label{Fhat.Green.Bound.lem} Let $f$ be a pdf supported on $[-\frac14,\frac14]$. For
$j\not=0$, $$|\hat{f}(j)|^2 \leq \frac{\ffi}{\pi} \,\sin\left(\frac{\pi}{\ffi}\right).$$
\end{lem}

\begin{proof}
Let $f_1:\T\to\R$ be defined by $f_1(x):=f(x-x_0)$, with $x_0$ chosen so that
$\hat{f_1}(j)$ is real and positive (clearly $\hat{f_1}(j)=|\hat{f}(j)|$ and
$\ffi=\|f_1\ast f_1\|_\infty$). Set $h(x)$ to be the symmetric decreasing rearrangement
of $f_1\ast f_1$, and $\overline{h}(x):=\ffi I(x)$, where $I(x)$ is the indicator
function of $[-\frac{1}{2 \ffi},\frac{1}{2 \ffi}]$. We have
 \begin{align*}
 |\hat{f}(j)|^2 &=      \hat{f_1}(j)^2
                =      \widehat{f_1\ast f_1}(j)
                =      \int f_1\ast f_1(x) \cos(2\pi j x)\,dx
                \leq    \int h(x) \cos(2\pi x)\,dx
 \end{align*}
by the inequality~\eqref{1988.Hardy.Littlewood.Polya.ineq}. We now apply Lemma~\ref{Focus.lem} to
find
 \begin{multline*}
 |\hat{f}(j)|^2
    \leq    \int \overline{h}(x)\cos(2\pi x)\,dx \\
    =       \int_{-1/(2\ffi)}^{1/(2\ffi)} \ffi \cos(2\pi x)\,dx
    =       \frac{\ffi}{\pi}\sin\left(\frac{\pi}{\ffi}\right).
 \end{multline*}
\end{proof}

With this technical result in hand, we can finally establish the lower bound on $\De$ given in Theorem~\ref{Delta.Summary.thm}(ii).

\begin{prop}
\label{the.full.De.bound.prop} $\De \ge \Dconstant \e^2$ for all $0\le\e\le1$.
\end{prop}

This gist of the proof of Proposition~\ref{the.full.De.bound.prop} is that if $\ffi$ is small, then
$\Re\hat{f}(1)$ is small by Lemma~\ref{Fhat.Green.Bound.lem}, and so $\|\ff\|_2^2$ is not very
small by Corollary~\ref{Cf.Main.Bound.cor}, whence $\ffi$ is not small. If $\ffi <
\twotimesDconstant$, then we get a contradiction.

\begin{proof}
Let $f$ be a pdf supported on $[-\frac14,\frac14]$, and assume that
    \begin{equation}\label{Main.Bound.0.eq}
    \ffi<\twotimesDconstant.
    \end{equation}
Set $x_1:=\Re\hat{f}(1)$. Since $f$ is supported on
$[-\frac14,\frac14]$, we see that $x_1>0$. By Lemma~\ref{Fhat.Green.Bound.lem},
    \begin{equation}\label{Main.Bound.1.eq}
    0< x_1 < 0.4191447.
    \end{equation}
However, we already know from Corollary~\ref{Cf.Main.Bound.cor} that
    \begin{equation}\label{Main.Bound.3.eq}
    \ffi\geq \|\ff\|_2 \geq
    1 +2 x_1^4 + (1.53890149 - 2.26425375 x_1)^4.
    \end{equation}
Routine calculus shows that there are no simultaneous solutions to the
inequalities~(\ref{Main.Bound.0.eq}), (\ref{Main.Bound.1.eq}), and
(\ref{Main.Bound.3.eq}). Therefore $\ffi\geq\twotimesDconstant$, whence Lemma~\ref{Fourier.Connection.lem} implies that
$\De\ge \Dconstant \e^2$.
\end{proof}

%
%

\subsection{The Kernel Problem}\label{kernel.problem.section}

Let $\K$ be the class of functions $K\in L^2$ satisfying $K(x)\ge1$ on
$[-\frac14,\frac14]$. Proposition~\ref{Cf.First.Result.prop} suggests the problem
of computing
\begin{equation*}
\inf_{K\in\K} \|\hat{K}\|_p =
 \inf_{K\in\K}
    \left(\sum_{j=-\infty}^\infty |\hat{K}(j)|^p\right)^{1/p}.
\end{equation*}
In Proposition~\ref{Cf.First.Result.prop} the case $p=\frac43$ arose, but using the
Hausdorff-Young inequality in place of Parseval's identity we are led to consider $1< p
\leq \frac43$. Also, we assumed in Proposition~\ref{Cf.First.Result.prop} that $K$ was
continuous, but this assumption can be removed by taking the pointwise limit of
continuous functions.

As similar problems occur in~\cite{Cilleruelo.Ruzsa.Trujillo} and in~\cite{2001.Green}, we feel it is worthwhile
to detail the thoughts and experiments that led to the kernel functions chosen in Corollaries~\ref{First.cor}
and~\ref{Cf.Main.Bound.cor}.

Our first observation is that if $G\in\K$, then so is
$K(x):=\frac12(G(x)+G(-x))$, and since $|\hat{K}(j)|=|\Re\hat{G}(j)|\leq
|\hat{G}(j)|$ we know that $\|\hat{K}\|_p\le\|\hat{G}\|_p$. Thus, we may
restrict our attention to the {\em even} functions in $\K$.

We also observe that $|\hat{K}(j)|$ decays more rapidly if many derivatives of $K$ are
continuous. This suggests that we should restrict our attention to continuous $K$,
perhaps even to infinitely differentiable $K$. However, computations suggest that the
best functions $K$ are continuous but {\em not} differentiable at $x=\frac14$ (see in
particular Figure~\ref{optimal.K.p43.fig} below).

In the argument of Proposition~\ref{Cf.First.Result.prop} we used the
inequality $\int f \leq \int fK$, which is an equality if we take $K$ to be
equal to 1 on $[-\frac14,\frac14]$, instead of merely at least 1. In light of
this, we should not be surprised if the optimal functions in $\K$ are exactly
1 on $[-\frac14,\frac14]$. This is supported by our computations.

Finally, we note that if $K_i\in\K$, and $\alpha_i>0$ with $\sum_i \alpha_i = 1$, then
$\sum_i \alpha_i K_i(x)\in \K$ also. This is particularly useful with $K_1(x):=1$.
Specifically, given any $K_2\in\K$ with known $\|\hat{K_2}\|_p$ (we stipulate
$\|K_2\|_1=\hat K(0)\le1$ to avoid technicalities), we may easily compute the
$\alpha\in[0,1]$ for which $\|\hat{K}\|_p$ is minimized, where $K(x):=\alpha K_1(x) +
(1-\alpha)K_2(x)$. We have
\begin{equation}
    \|\hat K\|_p^p = (\alpha + (1-\alpha)\hat K_2(0) )^p + (1-\alpha)^p
    \lnorm{1}{p}{\hat K}^p = (1-(1-\alpha)M)^p + (1-\alpha)^p N,
    \label{hatKpp.expression}
\end{equation}
where we have set $M := 1-\hat K_2(0)$ and $N := \lnorm{1}{p}{\hat K}^p$. Taking the
derivative with respect to $\alpha$, we obtain
\begin{equation*}
p(1-\alpha)^{p-1} \bigg( M \Big( \frac1{1-\alpha}-M \Big)^{p-1} - N \bigg),
\end{equation*}
the only root of which is $\alpha = 1 - \frac{M^{q/p}}{M^q + N^{q/p}}$ (where
$\frac1p+\frac1q=1$). It is straightforward (albeit tedious) to check by substituting
$\alpha$ into the second derivative of the expression~(\ref{hatKpp.expression}) that
this value of $\alpha$ yields a local maximum for $\|\hat K\|_p^p$. The maximum value
attained is then calculated to equal $N \big( M^q + N^{q/p} \big)^{1-p}$, which is
easily computed from the known function $K_2$.

Notice that when $p=\frac43$ (so $q=4$), applying Proposition
\ref{Cf.First.Result.prop} with our optimal function $K$ yields
\begin{equation*}
    \|\ff\|_2^2 \ge \|\hat K\|_{4/3}^{-4}
    = \Big( N \big( M^4 + N^3 \big)^{-1/3} \Big)^{-3} =
    \frac{M^4 + N^3}{N^3} = 1 + \frac{(1-\hat K_2(0))^4}{\lnorm{1}{4/3}{\hat K}^4},
\end{equation*}
whereupon we recover the conclusion of Proposition~\ref{Main.Bound.prop} with $m=1$.

Haar wavelets provide a convenient basis for $L^2([-\frac12,\frac12])$. We have numerically optimized the coefficients in various spaces of potential kernel functions $K$ spanned by short sums of Haar wavelets to minimize $\|K\|_{4/3}$ within those spaces. The resulting functions are shown in Figure~\ref{optimal.K.p43.fig}. This picture justifies restricting our attention to continuous functions that are constant on $[-\frac14,\frac14]$, and also implies that the optimal kernels are non-differentiable at $\pm\frac14$, indeed that their derivatives become unbounded near these points.

\begin{figure}
\begin{center}
 \begin{picture}(384,360)
    \put(24,108){$K(x)$}
    \put(24,228){$K(x)$}
    \put(24,348){$K(x)$}
    \put(215,108){$K(x)$}
    \put(214,228){$K(x)$}
    \put(213,348){$K(x)$}
    \put(178,64){$x$}
    \put(178,184){$x$}
    \put(178,304){$x$}
    \put(372,64){$x$}
    \put(372,184){$x$}
    \put(372,304){$x$}
    \put(85,88){$N=31$}
    \put(86,208){$N=7$}
    \put(86,328){$N=1$}
    \put(275,88){$N=63$}
    \put(275,208){$N=15$}
    \put(276,328){$N=3$}
    \includegraphics{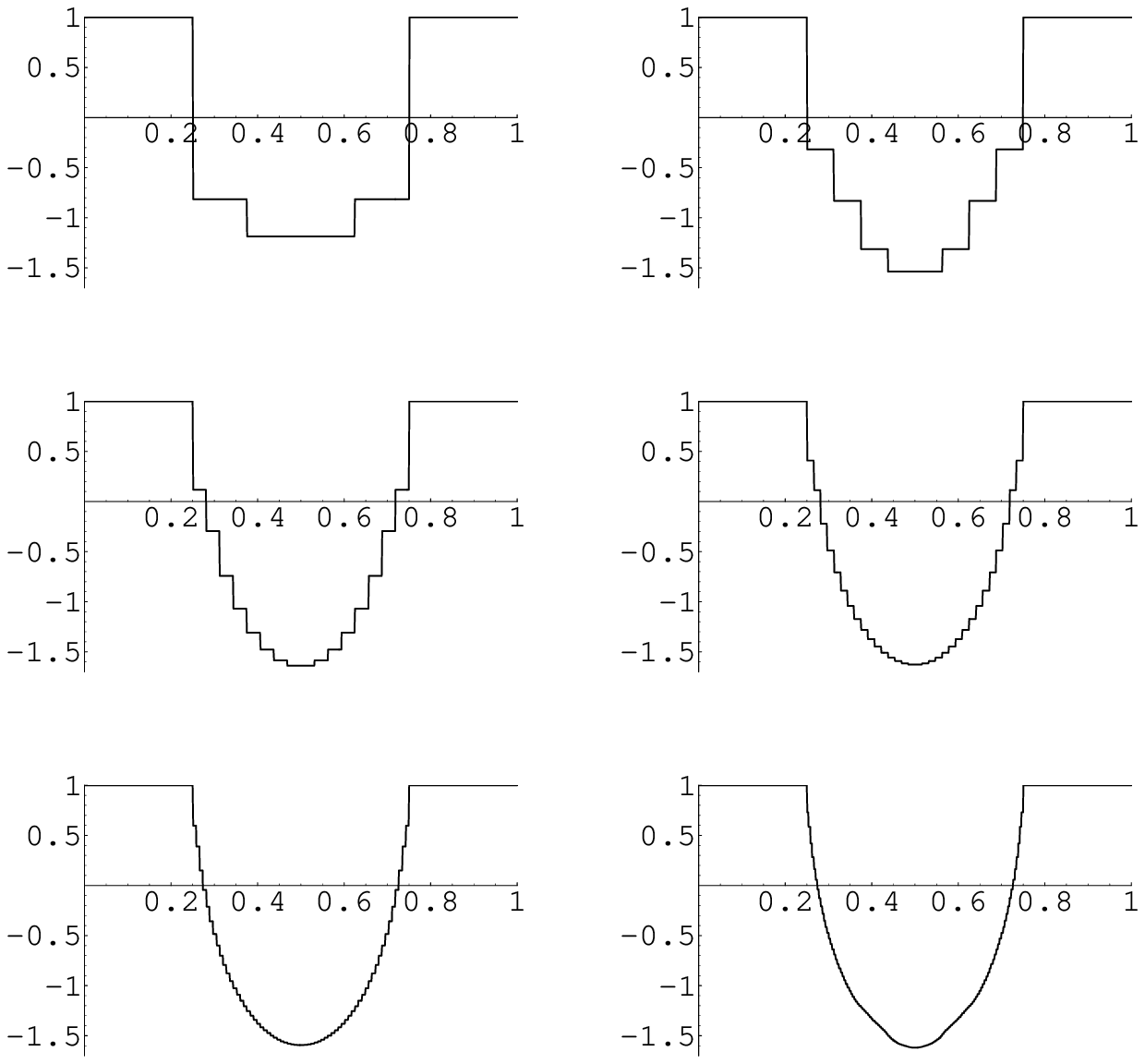}
 \end{picture}
 \begin{picture}(384,128)
    \put(25,122){$K(x)$}
    \put(186,102){$N=127$}
    \put(378,68){$x$}
    \put(0,0){\includegraphics{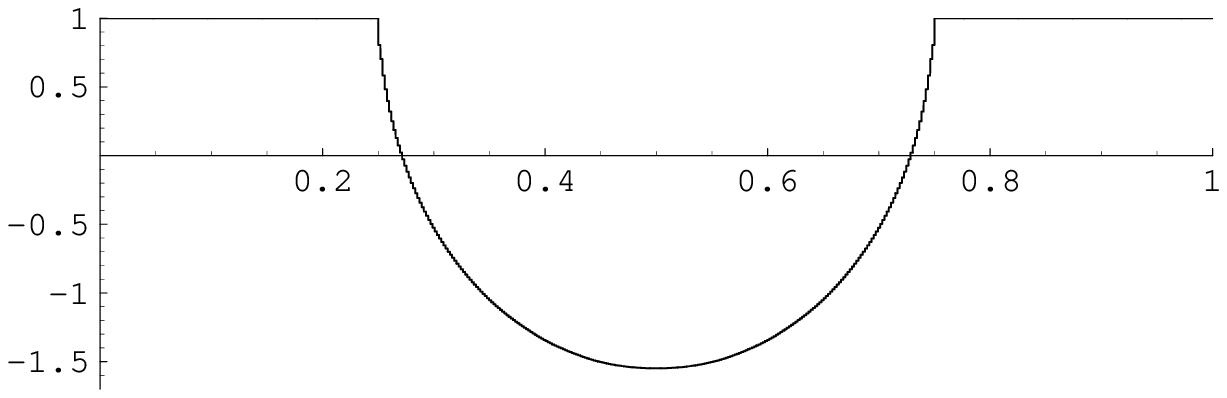}}
 \end{picture}
 \end{center}
    \caption{Optimal kernels generated by Haar wavelets\label{optimal.K.p43.fig}}
\end{figure}

For computational reasons, we further restrict attention to the class of continuous piecewise-linear even functions whose vertices all have abscissae with a given denominator. Let $\zeta(s,a)
:=\sum_{k=0}^\infty (k+a)^{-s}$ denote the Hurwitz zeta function. If ${\bf v}$ is a
vector, define $\Lambda_p({\bf v})$ to be the vector whose coordinates are the $p$th
powers of the absolute values of the corresponding coordinates of~${\bf v}$.

\begin{prop}
Let $T$ be a positive integer, $n$ a nonnegative integer, and $p\ge1$ a real number. For each integer $0\le t\le
T$, define $x_t:=\frac14+\frac{t}{4T}$, and let $y_t$ be an arbitrary real number, except that $y_0=1$. Let
$K(x)$ be the even function on $\T$ that is linear on $[0,\frac14]$ and on each of the intervals $[x_{t-1},x_t]$
$(1\leq t \leq T)$, satisfying $K(0)=1$ and $K(x_t)=y_t$ $(0\leq t \leq T)$. Then
    $$
    \lnorm{n}{p}{\hat{K}} =( 2 \Lambda_p(\vecd\A) \cdot\vecz )^{1/p},
    $$
where $\vecd$ is the $T$-dimensional vector $\vecd=(y_1-y_0, y_2-y_1, \dots,
y_T-y_{T-1})$, $\A$ is the $T\times4T$ matrix whose $(t,k)$-th component is
    $$
    \A_{tk} = \cos(2\pi(n+k-1)x_t)-\cos(2\pi(n+k-1)x_{t-1}),
    $$
and $\vecz$ is the $4T$-dimensional vector
    $$
    \vecz = (8T\pi^2)^{-p} \big( \zeta(2p,\tfrac{j}{4T}),
    \zeta(2p,\tfrac{j+1}{4T}), \dots, \zeta(2p,\tfrac{j+4T-1}{4T}) \big).
    $$
\label{Piecewise.Linear.prop}
\end{prop}

\begin{proof}
Note that
 \begin{align*}
   \hat{K}(-j)  =   \hat{K}(j)
                &=   \int_{-1/2}^{1/2} K(u)\cos(2\pi j u)\,du  \\
                &=   2\int_{0}^{1/4}\cos(2\pi j u)\,du+
                    2\sum_{t=1}^T\int_{x_{t-1}}^{x_t}(m_t u+b_t)\cos(2\pi j u)\,du,\\
 \end{align*}
where $m_t$ and $b_t$ are the slope and $y$-intercept of the line going
through $(x_{t-1},y_{t-1})$ and $(x_t,y_t)$. If we define $C(j):=\tfrac{\pi^2
j^2}{2T}\hat{K}(j)$, then integrating by parts we have
\begin{multline}
    C(j)= \left( \frac{\pi^2 j^2}{T}
    \left(\left. \frac{1}{2\pi j} \sin\left(2\pi j u\right)\right|_{0}^{1/4} \right)
    + \frac{\pi^2 j^2}{T}\sum_{t=1}^T \left( \left.  \frac{m_t u+b_t}{2\pi
    j}\sin(2\pi j u)\right|_{x_{t-1}}^{x_t}\right) \right) \\
    + \left( \frac{\pi^2 j^2}{T} \sum_{t=1}^T \left.\frac{m_t}{(2\pi j)^2}
    \cos(2\pi j u)\right|_{x_{t-1}}^{x_t} \right).
    \label{first.term.O1}
\end{multline}
The first term of this expression is
\begin{multline*}
    \frac{\pi j}{2T} \bigg( \sin\left( {\pi \tfrac j2} \right) + \sum_{t=1}^T
    \big(  (m_t x_t+b_t)\sin(2\pi j x_t) - (m_t x_{t-1}+b_t)\sin(2\pi
    j x_{t-1}) \big) \bigg) \\
    = \frac{\pi j}{2T} \bigg( \sin\left( {\pi \tfrac j2} \right) + \sum_{t=1}^T (m_t
    x_t+b_t)\sin(2\pi j x_t) - \sum_{t=0}^{T-1} (m_{t+1} x_t+b_{t+1})\sin(2\pi j x_t) \bigg).
\end{multline*}
Since $m_{t+1} x_t+b_{t+1} = y_t = m_t x_t+b_t$ and $x_0=\frac14$,
$x_T=\frac12$, this entire expression is a telescoping sum whose value is zero.
Eq.~(\ref{first.term.O1}) thus becomes
\begin{align}
    C(j) &= \frac{\pi^2 j^2}{T} \sum_{t=1}^T \left.\frac{m_t}{(2\pi j)^2} \cos(2\pi j
    u)\right|_{x_{t-1}}^{x_t} \notag \\
            &=   \sum_{t=1}^T    \left(y_t-y_{t-1}\right)
                                \left(\cos(2\pi j x_t)-\cos(2\pi j x_{t-1})\right)
                                \label{How.to.compute.C}
\end{align}
using $m_t=\frac{y_t-y_{t-1}}{x_t-x_{t-1}} = 4T(y_t-y_{t-1})$. Each $x_t$ is rational
and can be written with denominator $4T$, so we see that the sequence of normalized
Fourier coefficients $C(j)$ is periodic with period $4T$.

We proceed to compute $\lnorm{n}{p}{\hat{K}}$ with $n$ positive and
$p\geq1$.
\begin{align*}
    \left(\lnorm{n}{p}{\hat{K}}\right)^p
         &=  \sum_{|j|\geq n} |\hat{K}(j)|^p
         =  2\sum_{j=n}^\infty |\hat{K}(j)|^p
         =  2\sum_{j=n}^\infty \left|C(j)\frac{2T}{\pi^2 j^2}\right|^p \\
        &=  2\left(\frac{2T}{\pi^2}\right)^p
                    \sum_{j=n}^\infty \frac{|C(j)|^p}{j^{2p}}.
\end{align*}
Because of the periodicity of $C(j)$, we may write this as
\begin{align}
\left(\lnorm{n}{p}{\hat{K}}\right)^p
        &=  2\left(\frac{2T}{\pi^2}\right)^p
                    \left(\sum_{j=n}^{n+4T-1} |C(j)|^p
                            \sum_{r=0}^\infty (4Tr+j)^{-2p}\right) \notag \\
        &=  2\left(\frac{2T}{(4T)^2\pi^2}\right)^p
                    \left( \sum_{j=n}^{n+4T-1} |C(j)|^p
                            \zeta\left(2p,\tfrac{j}{4T}\right)
                            \right),\label{piecewise.linear.eq}
\end{align}
which concludes the proof.
\end{proof}

Proposition~\ref{Piecewise.Linear.prop} is useful in two ways. The first is that only
$\vecd$ depends on the chosen values $y_t$. That is, the vector $\vecz$ and the matrix
$\A$ may be precomputed (assuming $T$ is reasonably small), enabling us to compute
$\lnorm{n}{p}{\hat{K}}$ quickly enough as a function of $\vecd$ to numerically optimize
the $y_t$. The second use is through Eq.~(\ref{piecewise.linear.eq}). For a given $K$,
we set $y_t=K(x_t)$, whereupon $C(j)$ is computed for each $j$ using the formula in
Eq.~\eqref{How.to.compute.C}. Thus we can use Eq.~(\ref{piecewise.linear.eq}) to compute
$\lnorm{n}{p}{\hat{K_1}}$ with arbitrary accuracy, where $K_1$ is almost equal to $K$.
We have found that with $T=10000$ one can generally compute $\lnorm{n}{p}{\hat{K_1}}$
quickly.

In performing these numerical optimizations, we have found that ``good'' kernels $K(x)\in \K$ have a very
negative slope at $x=\frac14 ^+$. See Figure \ref{optimal.K.p43.fig}, for example, where the ``$N=2^j-1$''
picture denotes the step function $K$ for which $\|\hat K\|_{4/3}$ is minimal among all step functions whose
discontinuities all lie within the set $\frac1{2^{j+2}}\Z$.

Viewing graphs of these numerically optimized kernels suggests
that functions of the form
    $$
    K_{d_1,d_2}(x)=
    \begin{cases}
        1 & |x|\leq\tfrac14, \\
        1 - (1 - (4(\tfrac12 - x))^{d_1})^{d_2} & \tfrac14 < |x|\leq \tfrac12,
    \end{cases}
    $$
which have slope $-\infty$ at $x=\frac14^+$, may be very good. (Note that the graph of
$K_{2,1/2}(x)$ between $\frac14$ and $\frac34$ is the lower half of an ellipse.) More
good candidates are functions of the form
    $$
    K_{e_1,e_2,e_3}(x) =
    \begin{cases}
        1
        & |x|\leq\tfrac14, \\
        \left(\frac2\pi \tan^{-1}\left(\frac{(1-2x)^{e_1}}{(4x-1)^{e_2}}\right)\right)^{e_3}
        & \tfrac14 < |x|\leq \tfrac12,
    \end{cases}
    $$
where $e_1,e_2,$ and $e_3$ are positive. We have used a function of the form $K_{d_1,d_2}$ in the proof
of Corollary~\ref{Cf.Main.Bound.cor} and a function of the form $K_{e_1,e_2,e_3}$ in the proof of
Corollary~\ref{First.cor}.

\subsection{A Lower Bound for $\De$ around $\e=\frac12$}
\label{delta.half.section}

We begin with a fundamental relationship between $\Re \hat{f}(1)$ and $\Re \hat{f}(2)$.

\begin{lem}
    \label{Fhat.F1.F2.Domain.lem}
Let $f$ be a pdf supported on $[-\tfrac14,\tfrac14]$. Then
    $$2\big(\Re\hat{f}(1)\big)^2-1 \leq \Re\hat{f}(2) \leq 2(\Re\hat{f}(1))-1.$$
\end{lem}

\begin{proof}
To prove the first inequality, set $L_b(x)=b \cos(2\pi x)-\cos(4\pi x)$ (with $b\ge 0$) and observe
that for $-\tfrac14\leq x \leq\tfrac14$, we have $L_b(x)\leq 1+\frac{b^2}{8}$. Thus
    $$
    1+\tfrac{b^2}{8} \geq \int f(x)L_b(x)\,dx=\sum_{j=-2}^2
    \hat{f}(j)\hat{L_b}(-j)=b\Re\hat{f}(1)-\Re\hat{f}(2).
    $$
Rearranging, we arrive at $\Re\hat{f}(2) \geq b (\Re\hat{f}(1))-1-\frac{b^2}{8}$.
Setting $b=4\Re\hat f(1)$, we find that $\Re\hat{f}(2) \geq 2(\Re\hat{f}(1))^2-1$.

As for the second inequality, since $L_2(x):=2\cos(2\pi x)-\cos(4\pi x)$ is at least 1 for $-\tfrac14\leq x \leq\tfrac14$, we have
    $$1 \leq \int f(x)L(x)\,dx=\sum_{j=-2}^2
    \hat{f}(j)\hat{L}(-j)=2(\Re\hat{f}(1))-\Re\hat{f}(2).$$
Rearranging, we arrive at $\Re\hat{f}(2) \leq 2(\Re\hat{f}(1))-1$.
\end{proof}

From the inequality $\Re \hat{f}(2) \leq 2\Re\hat{f}(1)-1$
(Lemma~\ref{Fhat.F1.F2.Domain.lem}) one easily computes that $\max\{|\hat{f}(1)|,
|\hat{f}(2)|\} \geq \frac13$, and with Lemma~\ref{Fhat.Green.Bound.lem} this gives
    $$\frac19\leq \frac{\ffi}{\pi}\sin\left(\frac{\pi}{\ffi}\right).$$
This yields $\ffi \geq 1.11$, a non-trivial bound. If one assumes that $f$ is an nif
supported on a subset of $[-\frac14,\frac14]$ with large measure, then one can do much
better than Lemma~\ref{Fhat.F1.F2.Domain.lem}. The following proposition establishes the
lower bound on $\De$ given in Theorem \ref{Delta.Summary.thm}(iii).

\begin{prop}\label{Delta.one.half.thm}
Let $f$ be an nif supported on a subset of $[-\frac14,\frac14]$ with measure $\e/2$.
Then
    $$\ffi \geq 1.1092 + 0.176158\,\e $$
and consequently
    $$ \De \geq 0.5546\e^2 + 0.088079\e^3.$$
In particular, $\Delta(\frac12) \ge 0.14966$.
\end{prop}

\begin{proof}
For $\e\ge \frac58$, this proposition is weaker than
Lemma~\ref{Delta.Trivial.Bounds.lem}, and for $\e\le \frac38$ it is weaker than
Proposition~\ref{the.full.De.bound.prop}, so we restrict our attention to
$\frac38<\e<\frac58$.

Let $b>-1$ be a parameter and set $L_b(x):=\cos(4\pi x)-b \cos(2\pi x)$. If we define
$F:=\max\{\Re\hat{f}(1),-\Re\hat{f}(2)\}$, then
    $$
    \int f(x) L_b(x)\,dx = \Re\hat{f}(2)-b\Re\hat{f}(1) \geq -(b+1)F
    $$
on the one hand, and
    $$
    \int f(x) L_b(x)\,dx \leq \int \sdr{f}(x)\sdr{L_b}(x)\,dx
        =\int_{-\e/4}^{\e/4} \tfrac 2\e \,\sdr{L_b}(x)\,dx
    $$
on the other, where $\sdr{L_b}(x)$ is the symmetric decreasing rearrangement of $L_b(x)$
on the interval $[-\frac14,\frac14]$. Thus
    $$
    F \geq \frac{-1}{b+1} \frac 2\e \int_{-\e/4}^{\e/4} \sdr{L_b}(x)\,dx.
    $$
The right-hand side may be computed explicitly as a function of $\e$ and $b$ and then
the value of $b$ chosen in terms of $\e$ to maximize the resulting expression. One finds
that for $\e<\frac58$, the optimal choice of $b$ lies in the interval $2<b<4$, and the
resulting lower bound for $F$ is
    $$
    F \geq \frac{3\cos (\frac{\pi \e }{4}) +
    \sin (\frac{\pi \e }{4}) -
    {\sqrt{3 + 4\cos (\frac{\pi \e }{2}) +
        2\cos (\pi \e ) -
        \sin (\frac{\pi \e }{2})}}}{\pi
     \e \cos (\frac{\pi \e }{4}) +
    \pi \e \sin (\frac{\pi \e }{4})}.
    $$
From Lemma~\ref{Fhat.Green.Bound.lem} we know that $F^2\leq
\frac{\ffi}{\pi}\sin\left(\frac{\pi}{\ffi}\right) $. We compare these bounds on $F$ to
conclude the proof. Specifically,
    \begin{equation}\label{F2.lower.bound}
    F^2 \leq \frac{\ffi}{\pi}\sin\left(\frac{\pi}{\ffi}\right)
        \leq \frac{3}{5\pi } + \frac{\left( 6 + 5{\sqrt{3}}\pi  \right)
     \left( \ffi-\frac65 \right) }{12\pi },
    \end{equation}
where the expression on the right-hand side of this equation is from the Taylor
expansion of $\frac x\pi\sin(\frac\pi x)$ at $x_0=\frac65$, and
    \begin{align}
    F^2 &\geq \left(\frac{3\cos (\frac{\pi \e }{4}) +
         \sin (\frac{\pi \e }{4}) -
        {\sqrt{3 + 4\cos (\frac{\pi \e }{2}) +
        2\cos (\pi \e ) -
        \sin (\frac{\pi \e }{2})}}}{\pi
         \e \cos (\frac{\pi \e }{4}) +
        \pi \e \sin (\frac{\pi \e }{4})}\right)^2
            \notag \\
     &\geq \frac{-8\left( -3 - \sqrt{2} + \sqrt{3} + \sqrt{6} \right) }{ \pi^2}
            \notag \\
     &\qquad + \frac{\left( 96\left( -3 - {\sqrt{2}} + {\sqrt{3}} + {\sqrt{6}} \right)  -
       4\left( 9{\sqrt{2}} - 10{\sqrt{3}} + {\sqrt{6}} \right) \pi  \right)
     \left( \e-\frac12 \right) }{3{\pi }^2}, \label{F2.upper.bound}
     \end{align}
where the expression on the right-hand side is from the Taylor expansion of the middle
expression at $\e_0=\frac 12$. Comparing Eqs.~(\ref{F2.lower.bound})
and~(\ref{F2.upper.bound}) gives a lower bound on $\ffi$, say $\ffi\geq c_1 + c_2 \e$
with certain constants $c_1,c_2$. It is easily checked that $c_1> 1.1092 $ and
$c_2>0.176158$, concluding the proof of the first asserted inequality. The second
inequality then follows from Lemma \ref{Fourier.Connection.lem}.
\end{proof}

\section{Upper Bounds for $\De$}

\subsection{Inequalities Relating $\De$ and $R(g,n)$}

A symmetric set consists of pairs $(x,y)$ all with a fixed midpoint
$c=\tfrac{x+y}{2}$. If there are few pairs in $E\times E$ with a given sum
$2c$, then there will be no large symmetric subset of $E$ with center $c$. We
take advantage of the constructions of large integer sets whose pairwise sums
repeat at most $g$ times to construct large real subsets of $[0,1)$ with no
large symmetric subsets. More precisely, a set $S$ of integers is called a $\Bg$ set if for any given $m$ there are at most $g$ ordered pairs $(s_1,s_2)\in S \times S$ with $s_1+s_2=m$. (In the case $g=2$, these are better known as Sidon sets.) Define
\begin{equation}
R(g,n):= \max \big\{ |S| \colon S\subseteq \{1,2,\dots,n\},\, \text{$S$ is a $\Bg$ set} \big\}.
\label{Rdef}
\end{equation}

\begin{prop}
For any integers $n\ge g\ge1$, we have $\Delta(\frac{R(g,n)}n) \leq
\frac{g}n.$ \label{Delta.from.R.Connection.prop}
\end{prop}

\begin{proof}
Let $S\subseteq\{1,2,\dots,n\}$ be a $B^\ast[g]$ set with $|S|=R(g,n)$. Define
\begin{equation}
A(S) := \bigcup_{s\in S} \bigg[ \frac{s-1}{n},\frac{s}{n} \bigg),
\label{ASdef}
\end{equation}
and note that $A(S)\subseteq[0,1]$ and that the measure of $A(S)$ is exactly $R(g,n)/n$. Thus it
suffices to show that the largest symmetric subset of $A(S)$ has measure at most $\frac gn$.

Notice that the set $A(S)$ is a finite union of intervals, and so the
function $\lambda\big(A(S)\cap(2c-A(S))\big)$, which gives the measure of the largest
symmetric subset of $A(S)$ with center $c$, is piecewise linear. (Figure
\ref{AS.and.convolution} contains a typical example of the set $A(S)$ portrayed in dark
gray below the $c$-axis, together with the function $\lambda\big(A(S)\cap(2c-A(S))\big)$
shown as the upper boundary of the light gray region above the $c$-axis, for
$S=\{1,2,3,5,8,13\}$.) Without loss of generality, therefore, we may restrict our
attention to those symmetric subsets of $A(S)$ whose center $c$ is the midpoint of
endpoints of any two intervals $\big( \frac{s-1}{n},\frac{s}{n} \big)$. In other words,
we may assume that $2nc\in\Z$.

\begin{figure}[ht]
\begin{center}
    \begin{picture}(384,178)
        \put(0,163){$\lambda(A(S)\cap(2c-A(S)))$}
        \put(370,20){$c$}
        \put(30,60){$\tfrac{1}{13}$}
        \put(30,99.5){$\tfrac{2}{13}$}
        \put(30,139.5){$\tfrac{3}{13}$}
        \put(10,0){\includegraphics{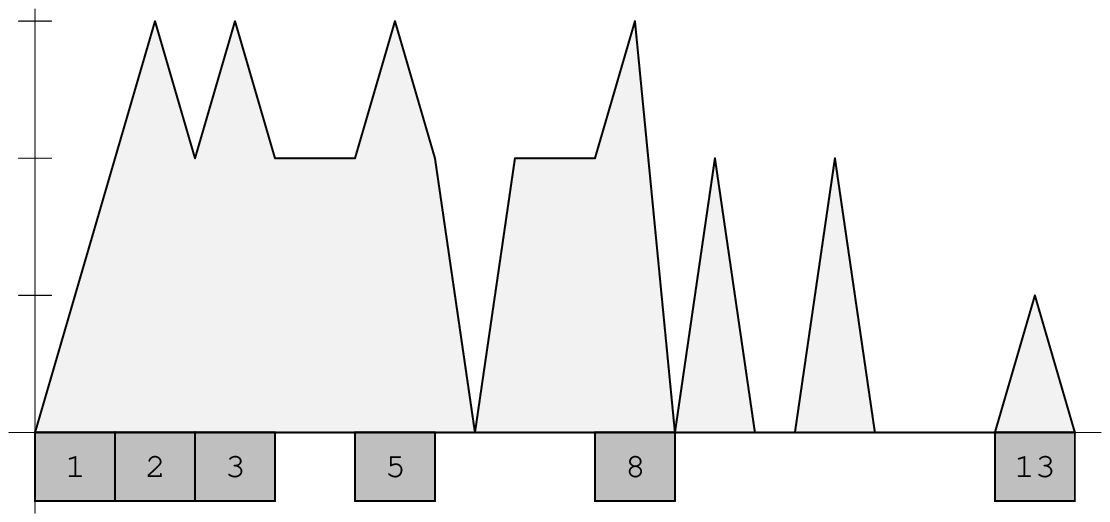}}
    \end{picture}
\end{center}
    \caption{$A(S)$, and the function $\lambda(A(S)\cap(2c-A(S)))$,
        with $S=\{1,2,3,5,8,13\}$\label{AS.and.convolution}}
\end{figure}

Suppose $u$ and $v$ are elements of $A(S)$ such that $\frac{u+v}2=c$. Write
$u=\frac{s_1}n - \frac1{2n} + x$ and $v=\frac{s_2}n - \frac1{2n} + y$ for integers
$s_1,s_2\in S$ and real numbers $x,y$ satisfying $|x|,|y|<\frac1{2n}$. (We may ignore
the possibility that $nu$ or $nv$ is an integer, since this is a measure-zero event for
any fixed $c$.) Then $2nc = n(u+v) = s_1+s_2-1+n(x+y)$, and since $2nc$, $s_1$, and
$s_2$ are all integers, we see that $n(x+y)$ is also an integer. But $|n(x+y)| < 1$, so
$x+y=0$ and $s_1+s_2 = 2nc+1$.

Since $S$ is a $\Bg$ set, there are at most $g$ solutions $(s_1,s_2)$ to the equation $s_1+s_2 = 2nc+1$. If it
happens that $s_1=s_2$, the interval $\big( \frac{s_1-1}{n},\frac{s_1}{n} \big)$ (a set of measure $\frac1n$) is
contributed to the symmetric subset with center $c$. Otherwise, the set $\big( \frac{s_1-1}{n},\frac{s_1}{n}
\big) \cup \big( \frac{s_2-1}{n},\frac{s_2}{n} \big)$ (a set of measure $\frac2n$) is contributed to the
symmetric subset with center $c$, but this counts for the two solutions $(s_1,s_2)$ and $(s_2,s_1)$. In total,
then, the largest symmetric subset having center $c$ has measure at most $\frac gn$. This establishes the
proposition.
\end{proof}

Using Proposition~\ref{Delta.from.R.Connection.prop}, we can translate lower bounds on
$\De$ into upper bounds on $R(g,n)$, as in Corollary~\ref{Delta.R.Basic.cor}.

\begin{cor}\label{Delta.R.Basic.cor}
If $\delta\le \inf_{0<\e<1} \De/\e^2$, then $R(g,n) \le \delta^{-1/2} \sqrt{gn}$ for all
$n\ge g \ge 1$.
\end{cor}

We remark that we may take $\delta = \Dconstant$ by Theorem~\ref{Delta.Summary.thm}(ii), and so
this corollary implies that $R(g,n) \le \Rconstant\sqrt{gn}$. This improves the previously-best
bound on $R(g,n)$ (given in \cite{2001.Green}) for $g\ge 30$ and $n$ large.

\begin{proof}
Combining the hypothesized lower bound $\Delta(\e)\ge \delta \e^2$ with
Proposition~\ref{Delta.from.R.Connection.prop}, we find that
    $$
    \delta\bigg( \frac{R(g,n)}n \bigg)^2 \le \Delta\bigg( \frac{R(g,n)}n \bigg)
    \le \frac gn
    $$
which is equivalent to $R(g,n)\le \delta^{-1/2}\sqrt{gn}$.
\end{proof}

We have been unable to prove or disprove that
 $$\lim_{g\to\infty} \lim_{n\to\infty} \frac{R(g,n)}{\sqrt{gn}} = \left(\inf_{0<\e<1}
 \frac{\De}{\e^2} \right)^{-1/2},$$
i.e., that Corollary~\ref{Delta.R.Basic.cor} is best possible as $g\to\infty$. At any rate, for small $g$ it is possible to
do better by taking advantage of the shape of the set $A(S)$ used in the proof of
Proposition~\ref{Delta.from.R.Connection.prop}. This is the subject of
the companion paper \cite{Martin.O'Bryant.b} of the authors.

Proposition~\ref{Delta.from.R.Connection.prop} provides a one-sided inequality linking
$\Delta(\e)$ and $R(g,n)$. It will also be useful for us to prove a theoretical result
showing that the problems of determining the asymptotics of the two functions are, in a
weak sense, equivalent. In particular, the following proposition implies
that the trivial lower bound $\Delta(\e)\ge\frac12\e^2$ and the trivial upper bound
$R(g,n)\le\sqrt{2gn}$ are actually equivalent. Further, any nontrivial lower bound on
$\De$ gives a nontrivial upper bound on $R(g,n)$, and vice versa.

\begin{prop}
$\Delta(\e) = \inf\{\frac gn\colon n\ge g\ge1,\, \frac{R(g,n)}n\ge\e\}$ for all
$0\le\e\le1$.
\label{Delta.as.infimum.prop}
\end{prop}

\begin{proof}
That $\Delta(\e)$ is bounded above by the right-hand side follows immediately from
Proposition~\ref{Delta.from.R.Connection.prop} and the fact that $\Delta$ is an increasing
function. For the complementary inequality, let $S\subseteq[0,1)$ with $\lambda(S)=\e$. Basic Lebesgue measure theory tells us that given any $\eta>0$, there
exists a finite union $T$ of open intervals such that $\lambda(S\oplus T) < \eta$, and it is easily
seen that $T$ can be chosen to meet the following criteria: $T\subseteq[0,1)$, the endpoints of the
finitely many intervals comprising $T$ are rational, and $\lambda(T)>\e$. Choosing a common
denominator $n$ for the endpoints of the intervals comprising $T$, we may write $T = \bigcup_{m\in
M} \big[ \frac {m-1}n, \frac mn \big)$ (up to a finite set of points) for some set of integers
$M\subseteq \{1,\dots,n\}$; most likely we have greatly increased the number of intervals
comprising $T$ by writing it in this manner, and $M$ contains many consecutive integers. Let $g$ be
the maximal number of solutions $(m_1,m_2)\in M\times M$ to $m_1+m_2=k$ as $k$ varies over all
integers, so that $M$ is a $B^\ast[g]$ set and thus $|M| \le R(g,n)$ by the definition of $R$. It
follows that $\e < \lambda(T) = |M|/n \le R(g,n)/n$. Now $T$ is exactly the set $A(M)$ as defined
in Eq.~\eqref{ASdef}; hence $D(T) = \frac gn$ as we saw in the proof of
Proposition~\ref{Delta.from.R.Connection.prop}. Therefore by Lemma~\ref{Diamond.lem},
\begin{equation*}
    D(S) \ge D(T) - 2\eta   = \tfrac gn - 2\eta  \ge
    \inf\{\tfrac gn\colon n\ge g \ge1,\,
                        \tfrac{R(g,n)}n\ge\e\} - 2\eta.
\end{equation*}
Taking the infimum over appropriate sets $S$ and noting that $\eta>0$ was
arbitrary, we derive the desired inequality $\Delta(\e) \ge \inf\{\frac
gn\colon n\ge g\ge1,\, \frac{R(g,n)}n\ge\e\}$.
\end{proof}

\subsection{Probabilistic Constructions of \Bgm Sets}
\label{circular.probabilistic.section}

We begin by considering a modular version of \Bg sets. A set $S$ is a \Bgm set if for any given $m$
there are at most $g$ ordered pairs $(s_1,s_2)\in S \times S$ with $s_1+s_2\equiv m\pmod n$
(equivalently, if the coefficients of the least-degree representative of $\left(\sum_{s\in S}
z^s\right)^2 \pmod{z^n-1}$ are bounded by $g$). For example, the set $\{0,1,2,4\}$ is a
$B^\ast[3]\pmod{7}$ set, and $\{0,1,3,7\}$ is a $B^\ast[2]\pmod{12}$ set. Note that $7+7\equiv 1+1
\pmod{12}$, so that $\{0,1,3,7\}$ is not a ``modular Sidon set'' as defined by some authors,
e.g.,~\cite{1980.2.Graham.Sloane} or \cite[Problem C10]{1994.Guy}.

Just as we defined $R(g,n)$ to be the largest possible cardinality of a \Bg set
contained in $[0,n)$, we define $C(g,n)$ to be the largest possible cardinality of a
\Bgm set. The mnemonic is ``R'' for the $\R$eal problem and ``C'' for the Circular
problem. We demonstrate the existence of large \Bgm sets via a probabilistic construction in this section, and we give a similar probabilistic construction of large \Bg sets in Section~\ref{probabilistic.bg.section}.

We rely upon the following two lemmas, which are quantitative statements of the Central Limit
Theorem.

\begin{lem}
Let $p_1,\dots,p_n$ be real numbers in the range $[0,1]$, and set
$p=(p_1+\dots+p_n)/n$. Define mutually independent random variables
$X_1,\dots,X_n$ such that $X_i$ takes the value $1-p_i$ with probability
$p_i$ and the value $-p_i$ with probability $1-p_i$ (so that the expectation
of each $X_i$ is zero), and define $X=X_1+\dots+X_n$. Then for any positive
number $a$,
\begin{equation*}
\Prob [ X>a ] < \exp \Big( \frac{-a^2}{2pn} + \frac{a^3}{2p^2n^2} \Big)
\quad\text{and}\quad
\Prob [ X<-a ] < \exp \Big( \frac{-a^2}{2pn} \Big).
\end{equation*}
\label{prob.method.lem}
\end{lem}

\begin{proof}
These assertions are Theorems A.11 and A.13 of~\cite{2000.Alon.Spencer}.
\end{proof}

\begin{lem}
Let $p_1,\dots,p_n$ be real numbers in the range $[0,1]$, and set
$E=p_1+\dots+p_n$. Define mutually independent random variables
$Y_1,\dots,Y_n$ such that $Y_i$ takes the value 1 with probability
$p_i$ and the value 0 with probability $1-p_i$, and define $Y=Y_1+\dots+Y_n$
(so that the expectation of $Y$ equals $E$). Then $\Prob [ Y>E+a ] < \exp
\big( \frac{-a^2}{3E} \big)$ for any real number
$0<a<E/3$, and $\Prob [ Y<E-a ] < \exp \big( \frac{-a^2}{2E} \big)$ for any
positive real number $a$.
\label{our.application.lem}
\end{lem}

\begin{proof}
This follows immediately from Lemma~\ref{prob.method.lem} upon defining $X_i = Y_i -
p_i$ for each $i$ and noting that $E=pn$ and that $\frac{a^3}{2E^2} < \frac{a^2}{6E}$
under the assumption $0<a<E/3$.
\end{proof}

We now give the probabilistic construction of large \Bgm sets.
We write that $f \gtrsim g$ (and $g\lesssim f$) if $\liminf f/g$ is at least 1.

\begin{prop}
For every $0<\e\le1$, there is a sequence of ordered pairs $(n_j,g_j)$ of positive
integers such that $\frac{C(g_j,n_j)}{n_j} \gtrsim \e$ and $\frac{g_j}{n_j}
\lesssim \e^2$.
\label{circle.probabilistic.prop}
\end{prop}

\begin{proof}
Let $n$ be an odd integer. We define a random subset $S$ of $\{1,\dots,n\}$ as follows:
for every $1\leq i \leq n$, let $Y_i$ be 1 with probability $\e$ and 0 with probability
$1-\e$ with the $Y_i$ mutually independent, and let $S := \{i \colon Y_i = 1\}$. We see that $|S|=\sum_{i=1}^n Y_i$ has
expectation $E=\e n$. Setting $a = \sqrt{\e n\log 4}$, Lemma~\ref{our.application.lem}
gives
$$\Prob \big[|S| < \e n - \sqrt{\e n\log 4}\, \big] < \frac12.$$

Now for any integer $k$, define the random variable
    \begin{align*}
    R_k &:= \# \{1\le c,d\le n \colon c+d\equiv k \pmod{n},\, Y_c=Y_d=1\}\\
        &= \sum_{c+d\equiv k \pmod{n}} Y_c Y_d,
    \end{align*}
so that $R_k$ is the number of representations of $k\pmod{n}$ as the sum of two elements
of $S$. Observe that $R_k$ is the sum of $n-1$ random variables taking the value 1 with
probability $\e^2$ and the value 0 otherwise, plus one random variable (corresponding to
$c\equiv d\equiv 2^{-1}k \pmod n$) taking the value 1 with probability $\e$ and the
value 0 otherwise. Therefore the expectation of $R_k$ is $E=(n-1)\e^2+\e$. Setting $a =
\sqrt{3((n-1)\e^2+\e)\log 2n}$, and noting that $a<E/3$ when $n$ is sufficiently large
in terms of $\e$, Lemma \ref{our.application.lem} gives
$$
\Prob \big[R_k > (n-1)\e^2+\e + \sqrt{3((n-1)\e^2+\e)\log 2n}\, \big] < \frac1{2n}
$$
for each $1\le k\le n$.

The random set $S$ is a $\Bgm$ set with $|S|>\e n-\sqrt{\e n \log 4}$ and $g \le E+a = (n-1)\e^2+\e
+ \sqrt{3((n-1)\e^2+\e)\log 2n}$ unless $|S|<\e n-\sqrt{\e n \log 4}$ or $R_1>E+a$ or $R_2>E+a$ or
$\ldots$ or $R_{2n}>E+a$. For any events $A_i$,
    \[\Prob [ A_1 \text{ or } A_2 \text{ or } \dots ]< \sum_{i} \Prob[A_i],\]
and consequently,
    \begin{multline*}
    \Prob \big[ \text{$|S|<\e n-\sqrt{\e n \log 4}$ or $R_1>E+a$ or $R_2>E+a$ or
    $\ldots$ or $R_{2n}>E+a$} \big] \\< \frac 12 + 2n\,\cdot \frac{1}{2n} < 1.
    \end{multline*}

Therefore, there exists a $\Bgm$ set $S\subseteq \{1,\dots,n\}$, with $g \le E+a=(n-1)\e^2+\e +
\sqrt{3((n-1)\e^2+\e)\log 2n} \lesssim \e^2n$, with $|S| \ge \e n - \sqrt{\e n\log 4} \gtrsim \e
n$. This establishes the proposition.
\end{proof}

Define $\Delta_\T(\e)$ to be the supremum of those real numbers $\delta$ such
that every subset of $\T$ with measure $\e$ has a subset with measure $\delta$
that is fixed by a reflection $t\mapsto c-t$. The function $\Delta_\T(\e)$
stands in relation to $C(g,n)$ as $\De$ stands to $R(g,n)$. However, it turns
out that $\Delta_\T$ is much easier to understand.

\begin{cor}
Every subset of $\T$ with measure $\e$ contains a symmetric subset
with measure $\e^2$, and this is best possible for every $\e$:
$$\Delta_\T(\e) = \e^2$$
for all $0\le\e\le1$.
\label{circle.answer.cor}
\end{cor}

\begin{proof}
In the proof of the trivial lower bound for $\De$ (Lemma
\ref{Delta.Trivial.Lower.Bound.lem}), we saw that every subset of $[0,1]$ with measure
$\e$ contains a symmetric subset with measure at least $\frac12\e^2$. The proof is
easily modified to show that every subset of $\T$ with measure $\e$ contains a symmetric
subset with measure $\e^2$. This shows that $\Delta_\T(\e) \ge \e^2$ for all $\e$. On
the other hand, the proof of Proposition~\ref{Delta.from.R.Connection.prop} is also
easily modified to show that $\Delta_\T\big( \frac{C(g,n)}n \big) \le \frac gn$, as is
the proof of Lemma~\ref{Lipschitz.Condition.lem} to show that $\Delta_\T$ is continuous.
Then, by virtue of Proposition~\ref{circle.probabilistic.prop} and the monotonicity of
$\Delta_\T$, we have $\Delta_\T(\e) \le \e^2$.
\end{proof}

\subsection{Probabilistic Constructions of \Bg Sets}
\label{probabilistic.bg.section}

We can use the probabilistic methods employed in Section
\ref{circular.probabilistic.section} to construct large \Bg sets in $\Z$. The proof is
more complicated because it is to our advantage to endow different integers with
different probabilities of belonging to our random set. Although all of the constants in
the proof could be made explicit, we are content with inequalities having error terms
involving big-O notation.

\begin{prop}
Let $\gamma\ge\pi$ be a real number and $n\ge\gamma$ be an integer. There exists a \Bg
set $S\subseteq\{1,\dots,n\}$, where $g = \gamma + O(\sqrt{\gamma\log n})$, with $|S|
\ge 2\sqrt{\frac{\gamma n}\pi} + O(\gamma + (\gamma n)^{1/4})$.
\label{probabilistic.bg.prop}
\end{prop}

\begin{proof}
Define mutually independent random variables $Y_k$, taking only the values 0 and 1, by
\begin{equation}
\Pr\{Y_k=1\} = p_k := \begin{cases}
1 & 1\le k < \frac\gamma\pi, \\
\sqrt{\!\frac\gamma{\pi k}} & \frac\gamma\pi \le k \le n, \\
0 & k > n.
\end{cases}
\label{schinzel.pk.def}
\end{equation}
(Notice that $p_k \le \sqrt{\!\frac\gamma{\pi k}}$ for all $k\ge1$.) These
random variables define a random subset $S = \{k\colon Y_k=1\}$ of the integers
from 1 to~$n$. We shall show that, with positive probability, $S$ is a large
$B^*[g]$ set with $g$ not much bigger than $\gamma$.

The expected size of $S$ is
\begin{align}
E_0 := \sum_{1\le j\le n} p_j &= \sum_{1\le j<\gamma/\pi} 1 +
\sum_{\gamma/\pi\le j\le n} \sqrt{\frac\gamma{\pi j}} \notag \\
&= \frac\gamma\pi + \int_{\gamma/\pi}^n \sqrt{\frac\gamma{\pi t}} \,dt + O(1)
= 2\sqrt{\frac{\gamma n}\pi} - \frac\gamma\pi + O(1).
\label{second.main.term.eq}
\end{align}
If we set $a_0 := \sqrt{2E_0\log3}$, then Lemma~\ref{our.application.lem} tells us that
    $$\Prob[|S| < E_0 - a_0] < \exp\big( \frac{-a_0^2}{2E_0} \big) = \frac13.$$

Now for any integer $k\in[\gamma,2n]$, let
\begin{equation*}
R_k := \sum_{1\le j\le n} Y_jY_{k-j} = 2\sum_{1\le j<k/2} Y_jY_{k-j} + Y_{k/2},
\end{equation*}
the number of representations of $k$ as $k=s_1+s_2$ with $s_1,s_2\in S$. (Here
we adopt the convention that $Y_{k/2}=p_{k/2}=0$ if $k$ is odd). Notice that
in this latter sum, $Y_{k/2}$ and the $Y_jY_{k-j}$ are mutually independent
random variables taking only values 0 and 1, with $\Pr[Y_jY_{k-j}=1] =
p_jp_{k-j}$. Thus the expectation of $R_k$ is
\begin{align}
E_k := 2\sum_{1\le j<k/2} p_jp_{k-j} + p_{k/2} &\le 2\sum_{1\le j<k/2}
\sqrt{\frac\gamma{\pi j}} \, \sqrt{\frac\gamma{\pi(k-j)}} +
\sqrt{\frac{\gamma}{\pi k/2}} \notag \\
&\le \frac{2\gamma}\pi \int_0^{k/2} \sqrt{\frac1{t(k-t)}} \,dt +
\sqrt{\frac{2\gamma}{\pi k}} = \gamma + \sqrt{\frac{2\gamma}{\pi k}} < \gamma+1
\label{pi.over.2.appearance}
\end{align}
using the inequalities $p_k \le \sqrt{\!\frac\gamma{\pi k}}$ and $k\ge\gamma$.

If we set $a = \sqrt{3(\gamma+1)\log3n}$, then Lemma~\ref{our.application.lem} tells us
that
\begin{equation*}
\Prob[R_k > \gamma + 1 + a] < \Prob[R_k > E_k + a] < \exp\Big(
\frac{-a^2}{3E_k} \Big) < \exp\Big( \frac{-a^2}{3(\gamma+1)} \Big) =
\frac1{3n}
\end{equation*}
for every $k$ in the range $\gamma\le k\le2n$. Note that $R_k\le\gamma$
trivially for $k$ in the range $1\le k\le\gamma$. Therefore, with probability
at least $1-\frac13-(2n-\gamma)\frac1{3n} = \frac\gamma{3n}>0$, the set $S$
has at least $E_0 - a_0 = 2\sqrt{\frac{\gamma n}\pi} + O(\gamma + (\gamma
n)^{1/4})$ elements and satisfies $R_k \le \gamma + 1 + a$ for all $1\le
k\le2n$. Setting $g:=\gamma+1+a = \gamma + O(\sqrt{\gamma\log n})$, we conclude
that any such set $S$ is a \Bg set. This establishes the proposition.
\end{proof}

Schinzel and Schmidt~\cite{2002.Schinzel.Schmidt} conjectured that among all pdfs supported on $[0,\frac12]$,
the function
\begin{equation*}
f(x) = \begin{cases}
\tfrac1{\sqrt{2x}} & x\in[0,\tfrac12], \\
0 & \text{otherwise}
\end{cases}
\end{equation*}
has the property that $\ffi$ is minimal. We have
\begin{equation*}
\ff(x) = \begin{cases}
\tfrac\pi2 & x\in[0,\tfrac12], \\
\tfrac\pi2 - 2\arctan \sqrt{2x-1} & x\in[\frac12,1], \\
0 & \text{otherwise}
\end{cases}
\end{equation*}
and so $\ffi = \frac\pi2$. We have adapted the function $f$ for our definition
(\ref{schinzel.pk.def}) of the probabilities $p_k$; the constant $\frac\pi2$ appears as
the value of the last integral in Eq.~(\ref{pi.over.2.appearance}). If Schinzel's
conjecture were false, then we could immediately incorporate any better function $f$
into the proof of Proposition~\ref{probabilistic.bg.prop} and improve the lower bound on
$|S|$.

\begin{thm}
For any $\delta>0$, we have $R(g,n)>\big( \frac2{\sqrt\pi} -\delta \big) \sqrt{gn}$ if both $\frac g{\log n}$ and $\frac ng$ are sufficiently large in terms of~$\delta$.
\label{g.growing.lower.bound.thm}
\end{thm}

\begin{proof}
In the proof of Proposition~\ref{probabilistic.bg.prop}, we saw that $\gamma\le g$ and
$g=\gamma+O(\sqrt{\gamma\log n})$; this implies that $\gamma=g+O(\sqrt{g\log n}) = g
\big( 1 + O\big( \sqrt{\frac{\log n}g} \big) \big)$. Therefore the size of the
constructed set $S$ was at least
\begin{align*}
2\sqrt{\tfrac{\gamma n}\pi} + O(\gamma + (\gamma n)^{1/4}) &= 2\sqrt{\tfrac
{gn}\pi \Big( 1 + O\Big( \sqrt{\tfrac{\log n}g} \,\Big) \Big)} +
O(g+(gn)^{1/4}) \\
&= 2\sqrt{\tfrac{gn}\pi} \Big( 1 + O \Big( \sqrt{\tfrac{\log n}g} +
\sqrt{\tfrac gn} \, \Big) \Big).
\end{align*}
This establishes the theorem.
\end{proof}

\subsection{Deriving the upper bounds}
\label{deriving.upper.bounds.sec}

In this section use the lower bounds on $R(g,n)$ established in Section~\ref{probabilistic.bg.section} into upper bounds for $\De$. Our first proposition verifies the statement of Theorem~\ref{Delta.Summary.thm}(i).

\begin{prop}
$\De=2\e-1$ for $\tfrac{11}{16}\leq \e \leq 1$.
\label{Line.2e-1.prop}
\end{prop}

\begin{proof}
We already proved in Lemma~\ref{Delta.Trivial.Bounds.lem} that $\De\ge2\e-1$ for all $0<\e\le1$.
Recall from Lemma~\ref{Lipschitz.Condition.lem} that the function $\Delta$ satisfies the Lipschitz condition $|\Delta(x)-\Delta(y)|\leq 2|x-y|$. Therefore to prove that $\De\le2\e-1$ for $\tfrac{11}{16}\leq \e \leq 1$, it suffices to prove simply that $\Delta\left(\frac{11}{16}\right) \le\frac38$.

For any positive integer $g$, it was shown by the authors \cite[Theorem 2(vi)]{Martin.O'Bryant.a} that
$$R(g,3g-\floor{g/3}+1) \ge g+2\floor{g/3}+\floor{g/6}.$$
We combine this with Proposition~\ref{Delta.from.R.Connection.prop} and the monotonicity of $\Delta$ to see that
\begin{equation*}
\frac{g}{3g-\floor{g/3}+1} \ge \Delta\bigg(
\frac{R(g,3g-\floor{g/3}+1)}{3g-\floor{g/3}+1} \bigg) \ge \Delta\bigg(
\frac{g+2\floor{g/3}+\floor{g/6}}{3g-\floor{g/3}+1} \bigg).
\end{equation*}
Since $\Delta$ is continuous by Lemma~\ref{Lipschitz.Condition.lem}, we may
take the limit of both sides as $g\to\infty$ to obtain
$\Delta\left(\frac{11}{16}\right) \le \frac38$ as desired.
\end{proof}

\noindent {\it Remark.} In light of the Lipschitz condition $|\Delta(x)-\Delta(y)|\leq 2|x-y|$, the lower bound $\De\ge2\e-1$ for all $0<\e\le1$ also follows easily from the trivial value $\Delta(1)=1$.

\begin{prop}
The function $\frac{\De}{\e^2}$ is increasing on $(0,1]$.
\label{De.over.e2.increasing.prop}
\end{prop}

\begin{proof}
The starting point of our proof is the inequality \cite[Theorem 2(v)]{Martin.O'Bryant.a}
$$R(g,x) C(h,y) \leq R(gh,xy).$$
With the monotonicity of $\De$ and Proposition~\ref{Delta.from.R.Connection.prop}, this gives
$$
\Delta \left( \frac{R(g,x)}{x} \frac{C(h,y)}{y} \right) \leq \Delta\left(
\frac{R(gh,x y)}{x y}  \right) \leq\frac{gh}{xy}.
$$
Choose $0<\e<\e_0$. Let $g_i,x_i$ be such that $\frac{R(g_i,x_i)}{x_i}\to \e_0$ and
$\frac{g_i}{x_i}\to \Delta(\e_0)$, which is possible by
Proposition~\ref{Delta.as.infimum.prop}. By Proposition
\ref{circle.probabilistic.prop}, we may choose sequences of integers $h_j$ and
$y_j$ such that $\frac{C(h_j,y_j)}{y_j} \gtrsim \frac\e{\e_0}$ and
$\frac{h_j}{y_j} \lesssim \big( \frac\e{\e_0} \big)^2$ as $j\to\infty$. This implies
$$
\frac{R(g_i,x_i)}{x_i} \frac{C(h_j,y_j)}{y_j} \gtrsim \e
 \quad \text{and}\quad
\frac{g_i}{x_i} \frac{h_j}{y_j} \lesssim \Delta(\e_0) \Big( \frac{\e}{\e_0}
\Big)^2,
$$
so that, again using the monotonicity and continuity of $\Delta$,
\begin{equation*}
\Delta(\e_0) \frac{\e^2}{\e_0^2} \gtrsim \frac{g_ih_j}{x_iy_j} \ge \Delta
\left( \frac{R(g_i,x_i)}{x_i} \frac{C(h_j,y_j)}{y_j} \right) \gtrsim \De
\end{equation*}
as $j\to\infty$. This shows that $\frac{\De}{\e^2} \le \frac{\Delta(\e_0)}{\e_0^2}$ as desired.
\end{proof}

We can immediately deduce two nice consequences of this proposition.

\begin{cor}
$\lim_{\e\to0^+} \frac{\De}{\e^2}$ exists.
\end{cor}

\begin{proof}
This follows from the fact that the function $\frac{\De}{\e^2}$ is
increasing and bounded below by $\frac12$ on $(0,1]$ by the trivial lower
bound (Lemma \ref{Delta.Trivial.Lower.Bound.lem}).
\end{proof}

\begin{cor}
$\De \leq \tfrac{96}{121}\e^2$ for $0\leq \e\leq \frac{11}{16}$.
\label{part.iii.cor}
\end{cor}

\begin{proof}
This follows from the value $\Delta\big( \frac{11}{16} \big) = \frac38$
calculated in Proposition~\ref{Line.2e-1.prop} and the fact that the function
$\frac{\De}{\e^2}$ is increasing.
\end{proof}

The corollary above proves part (iv) of Theorem~\ref{Delta.Summary.thm}, leaving only part (v) yet to be established. The following proposition finishes the proof of Theorem~\ref{Delta.Summary.thm}.

\begin{prop}
$\frac\De{\e^2} \le \frac\pi{(1+\sqrt{1-\e})^2}$ for all $0<\e\le1$.
\label{part.iv.prop}
\end{prop}

\begin{proof}
Define $\alpha:=1-\sqrt{1-\e}$, so that
$2\alpha-\alpha^2=\e$.
If we set $\gamma=\pi\alpha^2 n$ in the proof of Proposition
\ref{probabilistic.bg.prop}, then the sets constructed are \Bg sets with $g =
\pi\alpha^2n + O(\sqrt{n\log n})$ and have size at least
\begin{equation*}
E_0 - a_0 = 2\sqrt{\frac{\pi\alpha^2n^2}\pi} - \frac{\pi\alpha^2n}\pi + O(1 +
a_0) = (2\alpha-\alpha^2)n + O((\gamma n)^{1/4}) = \e n + O(\sqrt n)
\end{equation*}
from Eq.~(\ref{second.main.term.eq}).

Therefore, for these values of $g$ and $n$,
\begin{equation*}
\Delta \big( \frac{R(g,n)}n \big) \ge \Delta \big( \frac{\e n + O(\sqrt n)}n
\big) \to \De
\end{equation*}
as $n$ goes to infinity, by the continuity of $\Delta$. On the other hand, we see by
Proposition~\ref{Delta.from.R.Connection.prop} that
\begin{multline*}
\e^{-2}\Delta \big( \frac{R(g,n)}n \big) \le \frac g{\e^2n} =
\frac{\pi\alpha^2n + O(\sqrt{n\log n})}{\e^2 n} \\
= \frac{\pi \alpha^2}{(2\alpha-\alpha^2)^2} + O\Big( \sqrt{\frac{\log n}{\e^2
n}} \,\Big) =  \frac{\pi}{(2-\alpha)^2} + o(1) \to
\frac{\pi}{(1+\sqrt{1-\e})^2}
\end{multline*}
as $n$ goes to infinity. Combining these two inequalities yields
$\frac{\De}{\e^2} \le \frac{\pi}{(1+\sqrt{1-\e})^2}$ as desired.
\end{proof}

\section{Some Remaining Questions}

We group the problems in this section into three categories, although some problems do
not fit clearly into any of the categories and others fit into more than one.

\subsection{Properties of the Function $\De$}

The first open problem on the list must of course be the exact determination
of $\De$ for all values $0\le\e\le1$. In the course of our investigations, we
have come to believe the following assertion.

\begin{cnj}
$\De = \max\{2\e-1,\frac\pi4\e^2\}$ for all $0\le\e\le1$.
\label{bold.De.cnj}
\end{cnj}

Notice that the upper bounds given in Theorem~\ref{Delta.Summary.thm} are not too far
from this conjecture, the difference between the constants $\frac{96}{121} \doteq 0.79339$ and
$\frac\pi4 \doteq 0.78540$ in the middle range for $\e$ being the only discrepancy. In fact, we
believe it might be possible to prove that the expression in Conjecture
\ref{bold.De.cnj} is indeed an upper bound for $\De$ by a more refined application of
the probabilistic method employed in Section~\ref{probabilistic.bg.section}. The key
would be to show that the various events $R_k > \gamma + 1 + a$ are more or less
independent of one another (as it stands we have to assume the worst---that they are all
mutually exclusive---in obtaining our bound for the probability of obtaining a ``bad'' set).

There are some intermediate qualitative results about the function $\De$ that
might be easier to resolve. It seems likely that $\De$ is convex, for example,
but we have not been able to prove this. A first step towards clarifying the
nature of $\De$ might be to prove that
$$\frac{|\Delta(x)-\Delta(y)|}{|x-y|} \ll \max\{x,y\}.$$
Also, we would not be surprised to see accomplished an exact computation of
$\Delta(\tfrac12)$, but we have been unable to make this computation ourselves.
We do at least obtain $\Delta(\tfrac12) \ge 0.14966$ in Proposition~\ref{Delta.one.half.thm}. Note that Conjecture~\ref{bold.De.cnj} would imply that $\Delta(\tfrac12)=\frac\pi{16} \doteq 0.19635$.

We do not believe that there is always a set with measure $\e$ whose largest symmetric
subset has measure precisely $\De$. In fact, we do not believe that there is a set with measure
$\e_0:=\inf\{\e\colon \De = 2\e-1\}$ whose largest symmetric subset has measure
$\Delta(\e_0)$, but we do not even know the value of $\e_0$. In
Proposition~\ref{Line.2e-1.prop}, we showed that $\e_0 \leq \tfrac{11}{16}$, but this
was found by rather limited computations and is unlikely to be sharp. The quantity $\frac{96}{121}$ in Theorem
\ref{Delta.Summary.thm}(iv) is of the form $\frac{2\e_0-1}{\e_0^2}$, and thus any
improvement in the bound $\e_0\leq \frac{11}{16}$ would immediately result in
an improvement to Theorem~\ref{Delta.Summary.thm}(iv). We remark that Conjecture \ref{bold.De.cnj} implies
that $\e_0 = \frac2{2+\sqrt{4-\pi}} \doteq 0.68341$, which in turn would allow us to replace
the constant $\frac{96}{121} \doteq 0.79339$ in Theorem~\ref{Delta.Summary.thm}(iv) by
$\frac\pi4 \doteq 0.78540$.

\subsection{Artifacts of our Proof}\label{artifacts.section}

Let $\K$ be the class of functions $K\in L^2(\T)$ satisfying $K(x)\ge1$ on
$[-\frac14,\frac14]$. How small can we make $\|\hat{K}\|_p$ for $1\le p\le2$?
We are especially interested in $p=\frac43$, but a solution for any $p$ may be
enlightening.

To give some perspective to this problem, note that a trivial upper bound for
$\inf_{K\in\K} \|\hat{K}\|_p$ can be found by taking $K$ to be identically equal to 1,
which yields $\|\hat K\|_p=1$. One can find functions that improve upon this trivial
choice; for example, the function $K$ defined in Eq.~(\ref{easy.K.def}) is an example
where $\|\hat{K}\|_{4/3} \doteq 0.96585$. On the other hand, since the $\ell^p$-norm of a
sequence is a decreasing function of $p$, Parseval's identity immediately gives us the
lower bound $\|\hat K\|_p \ge \|\hat K\|_2 = \|K\|_2 \ge \big( \int_{-1/4}^{1/4} 1^2\,dt
\big)^{1/2} = \frac1{\sqrt2}  \doteq 0.70711$, and of course $\frac1{\sqrt2}$ is the exact minimum for
$p=2$.

We remark that Proposition~\ref{Cf.First.Result.prop} and the function $b(x)$ defined
after the proof of Corollary~\ref{First.cor} provide a stronger lower bound for $1\le p\le\frac43$. By direct
computation we have $1.14939> \|b\ast b\|_2^2$, and by
Proposition~\ref{Cf.First.Result.prop} we have $\|b\ast b\|_2^2 \geq
\|\hat{K}\|_{4/3}^{-4}$ for any $K\in\K$. Together these imply that
$\|\hat{K}\|_{p}\geq \|\hat{K}\|_{4/3}> 0.96579.$  In particular, for $p=\frac43$ we know the value of $\inf_{K\in\K} \|\hat K\|_{4/3}$ to within one part in ten thousand.
The problem of determining the actual
infimum for $1<p<2$ seems quite mysterious. We remark that Green~\cite{2001.Green}
considered the discrete version of a similar optimization problem, namely the
minimization of $\|\hat K\|_p$ over all pdfs $K$ supported on $[ {-\frac14},\frac14 ]$.

As mentioned at the end of Section~\ref{basic.argument.section}, we used the inequality $\| g \|_2^2 \leq
\|g\|_\infty \|g\|_1$ which is exact when $g$ takes on one non-zero value, i.e., when $g$ is an nif. We apply
this inequality when $g=\ff$ with $f$ supported on an interval of length $\frac12$, which usually looks very
different from an nif. In this circumstance, the inequality does not seem to be best possible, although the
corresponding inequality in the exponential sums approach of~\cite{Cilleruelo.Ruzsa.Trujillo} and in the
discrete Fourier approach of \cite{2001.Green} clearly is best possible. Specifically, we ask for a lower bound
on
 $$ \inf_{f:\R\mapsto \R_{\ge0}} \frac{\ffi \|\ff\|_1}{\| \ff \|_2^2}$$
that is strictly greater than 1. We know that this infimum is at most $\frac\pi{\log 16} \doteq 1.1331$, and in
fact Conjecture~\ref{infinity.two.cnj} would imply that the infimum is exactly $\frac\pi{\log 16}$.

\subsection{The Analogous Problem for Other Sets}

More generally, for any subset $E$ of an abelian group endowed with a measure, we can
define $\Delta_E(\e):=\inf\{D(A) \colon A \subseteq E,\, \lambda(A)=\e\}$, where $D(A)$
is defined in the same way as in Eq.~(\ref{Ddef}). For example, $\Delta_{[0,1]}(\e)$ is
the function $\De$ we have been considering throughout this paper, and $\Delta_\T(\e)$
was considered in Section \ref{circular.probabilistic.section}.

Most of the work in this paper generalizes easily from $E=[0,1]$ to
$E=[0,1]^d$. We have had difficulties, however, in finding good kernel
functions in higher dimensions. That is, we need functions $K(\bar{x})$ such
that
$$\sum_{\bar{j}\in\Z^d} \big| \hat{K} \big( \bar{j} \big) \big|^{4/3}$$
is as small as possible, while $K(\bar{x})\geq 1$ if all components of
$\bar{x}$ are less than $\frac14$ in absolute value. This restricts $K$ on
one-half of the space in 1 dimension, one-quarter of the space in 2 dimensions,
and only $2^{-d}$ of the space in $d$ dimensions. For this reason one might
expect that {\em better} kernels exist in higher dimensions, but the
computational difficulties have prevented us from finding them.

\bigskip{\small
{\it Acknowledgements.} The authors thank Heini Halberstam for thoughtful
readings of this man\-u\-script and helpful suggestions. The first author was
supported in part by grants from the Natural Sciences and Engineering Research Council.
The second author was supported by an NSF--Vigre Fellowship and grant DMS-0202460.}


\begin{thebibliography}{BVV00}


\bibitem[Alon and Spencer 2000]{2000.Alon.Spencer}
Noga Alon and Joel~H. Spencer.
\newblock {\em The probabilistic method}.
\newblock Wiley-Interscience [John Wiley \& Sons], New York, second edition,
  2000.
\newblock With an appendix on the life and work of Paul Erd\H os.

\bibitem[Beckner 1975]{Beckner}
W.~Beckner.
\newblock Inequalities in {F}ourier analysis.
\newblock {\em Ann. of Math.}, 102:159--182, 1975.



\bibitem[Banakh et al. 2000]{2000.Banakh.Verbitsky.Vorobets}
T.~Banakh, O.~Verbitsky, and Ya. Vorobets.
\newblock A {R}amsey treatment of symmetry.
\newblock {\em Electron. J. Combin.}, 7(1):Research Paper 52, 25 pp.
  (electronic), 2000.

\bibitem[Chung et al. 2000]{Chung.Erdos.Graham}
Fan Chung, Paul Erd{\H{o}}s, and Ronald Graham.
\newblock On sparse sets hitting linear forms.
\newblock {\em Number theory for the millennium, I}, 257--272, 2000.


\bibitem[Cilleruelo et al. 2002]{Cilleruelo.Ruzsa.Trujillo}
J.~Cilleruelo, I.~Ruzsa, and C.~Trujillo.
\newblock Upper and lower bounds for finite ${B}_h[g]$ sequences, $g>1$.
\newblock {\em J. Number Theory} 97, no. 1, 26--34, 2002.





\bibitem[Folland 1984]{1984.Folland}
Gerald~B. Folland.
\newblock {\em Real analysis}.
\newblock John Wiley \& Sons Inc., New York, 1984.
\newblock Modern techniques and their applications, A Wiley-Interscience
  Publication.

\bibitem[Green 2001]{2001.Green}
Ben Green.
\newblock The number of squares and ${B}_h[g]$ sets.
\newblock {\em Acta Arithmetica}, 100(4):365--390, 2001.

\bibitem[Graham and Sloane 1980]{1980.2.Graham.Sloane}
R.~L. Graham and N.~J.~A. Sloane.
\newblock On additive bases and harmonious graphs.
\newblock {\em SIAM J. Algebraic Discrete Methods}, 1(4):382--404, 1980.

\bibitem[Guy 1994]{1994.Guy}
Richard~K. Guy.
\newblock {\em Unsolved problems in number theory}.
\newblock Springer-Verlag, New York, second edition, 1994.
\newblock Unsolved Problems in Intuitive Mathematics, I.

\bibitem[Hardy et al. 1988]{1988.Hardy.Littlewood.Polya}
G.~H. Hardy, J.~E. Littlewood, and G.~P\'{o}lya.
\newblock {\em Inequalities}.
\newblock Cambridge University Press, Cambridge, 1988.
\newblock Reprint of the 1952 edition.



\bibitem[Martin and O'Bryant 2006]{Martin.O'Bryant.a}
Greg Martin and Kevin O'Bryant.
\newblock {Constructions of Generalized Sidon Sets}.
\newblock {\em J. Comb. Thy. Ser. A}, 113(4):591--607, 2006.

\bibitem[Martin and O'Bryant]{Martin.O'Bryant.b}
Greg Martin and Kevin O'Bryant.
\newblock {Upper Bounds for Generalized Sidon Sets}.
\newblock In preparation.


\bibitem[O'Bryant 2004]{2004.O'Bryant}
Kevin O'Bryant.
\newblock A Complete Annotated Bibliography of Work Related to Sidon Sequences.
\newblock {\em Elec. J. Combin.}, DS11, {\tt http://www.combinatorics.org/Surveys/}, 2004.




\bibitem[Schinzel and Schmidt 2002]{2002.Schinzel.Schmidt}
A. Schinzel and W. M. Schmidt.
\newblock Comparison of $L^1-$ and $L^\infty-$norms of squares of polynomials.
\newblock {\em Acta Arithmetica}, 104(3):283--296, 2002.

\bibitem[{\'S}wierczkowski 1958]{1958.Swierczkowski}
S.~{\'S}wierczkowski.
\newblock On the intersection of a linear set with the translation of its
  complement.
\newblock {\em Colloq. Math.}, 5:185--197, 1958.


\end{thebibliography}
\end{document}